\documentclass[12pt]{iopart}

\usepackage[utf8]{inputenc}
\usepackage{subfigure}
\usepackage{times}
\usepackage{epsfig,graphics,float,multirow} 
\usepackage{graphicx}
\usepackage{setstack}
\usepackage{iopams}
\usepackage[colorlinks=true, hyperfigures=false,
           urlcolor=blue, 
           pdfhighlight =/O]{hyperref}
\usepackage{float}

\newcommand{\R}{\mathbb{R}}

\newtheorem{Proposition}{Proposition}
\newtheorem{Lemma}{Lemma}
\newtheorem{Remark}{Remark}
\def\eps{\varepsilon}
\begin{document}
\title[Weakly coupled two slow- two fast systems]{Weakly coupled two slow- two fast systems, folded singularities and mixed mode oscillations}
\author{M. Krupa$^1$, B. Ambrosio, M.A. Aziz-Alaoui$^2$ }

\address{$1$ INRIA Paris-Rocquencourt Research Center, Domaine de Voluceau BP 105, 78153 Le Chesnay cedex, France}
\address{$2$ Normandie Univ, France;  ULH, LMAH, F-76600 Le Havre; FR CNRS 3335,  25 rue Philippe Lebon 76600 Le Havre, France.}
\ead{benjamin.ambrosio@univ-lehavre.fr}
\begin{abstract}
We study Mixed Mode Oscillations (MMOs) in systems of two weakly coupled slow/fast oscillators. We focus on the existence
and properties of a folded singularity called FSN II
that allows the emergence of MMOs in the presence of a suitable global return mechanism. 
As FSN II corresponds to a transcritical bifurcation for a desingularized reduced system, we prove that, under certain non-degeneracy conditions, such a transcritical bifurcation exists.
We then apply this result to the case of two coupled systems of FitzHugh-Nagumo type. This leads to a non trivial condition on the coupling that enables the existence of MMOs. 

\end{abstract}

\maketitle
\section{Introduction}
{\it Canard phenomenon}, known also as {\em canard explosion} is a transition from a small, Hopf type oscillation to a relaxation oscillation, occurring upon variation of a parameter. This transition has been first found in the context of  the van der Pol equation
\cite{eb-jlc-fd-md_81}, and soon after in numerous models of phenomena occurring in engineering and in chemical reactions 
\cite{mb_88}.
A common feature of all these models is the presence of time scale separation (one slow, one fast variable). 
A particular feature of canard explosion is that it takes place in a very small parameter interval. 
For the van der Pol system, where the
ratio of the timescales is given by a small parameter $\eps$, the width of this parameter interval can be estimated by
$O(\exp(-c/\eps)$, where $c>0$ is a fixed positive constant \footnote{Strictly speaking, if one defines canard explosion transition between small canards ans canard cycles with large head, then the width of the parameter interval where the transition takes place is given by $O(\exp(-c/\eps)$. }. The transition occurs through a sequence of {\it canard cycles}
whose striking feature is that they contain segments following closely a stable slow manifold and subsequently an unstable
slow manifold.  

The work on canard explosion led to investigations of slow/fast systems in three dimensions, with two slow and
one fast variables. In the context of these systems the concept of a {\em canard solution} 
or simply {\em canard} has been introduced,
as a solution passing from a stable to an unstable slow manifold \cite{eb_90,am-ps-hl-eg_98,ps-mw_01}.
Canards arise near so called folded singularities of which the best known is folded node  \cite{eb_90,ps-mw_01}
and \cite{mw_05}. Unlike in systems with one slow variable, canards occur robustly (in $O(1)$ parameter regions) in systems with two slow variables. Related
to canards are Mixed Mode Oscillations (MMOs), which are trajectories that combine small oscillations and large oscillations of relaxation type, both recurring in an alternating manner. Recently there has been a lot of interest in MMOs that arise due to the presence of canards, starting with the work of Milik, Szmolyan, Loeffelmann and Groeller \cite{am-ps-hl-eg_98}.  The small oscillations arise during the passage of the trajectories near a folded node or a related folded singularity. The dynamics near the folded singularity is transient, yet recurrent: the trajectories return to the neighborhood of the folded singularity by way of a global return mechanism \cite{mb-mk-mw_06}.

The setting of folded node combined with a global return mechanism, elucidated in \cite{mb-mk-mw_06}, led to the explanations of MMO dynamics found in applications \cite{rw_07, hr-mw-nk_08, be-mw_09, tv-rb-jt-mw_10}. A shortcoming of the folded node setting is the lack of connection to a Hopf bifurcation, which seems to play a prominent role in many MMOs.  This led to the interest in another, more degenerate, folded singularity, known as \textit{Folded Saddle Node of type II} (FSN II), originally introduced in \cite{am-ps-hl-eg_98} and recently analyzed in some detail by Krupa and Wechselberger \cite{mk-mw_10}. Guckenheimer \cite{jg_08} studied a very similar problem in the parameter regime yet closer to the Hopf bifurcation, calling it \textit{singular Hopf bifurcation}. The transition between the two settings was 
studied in \cite{rc-jr_11}. 
Another interesting singularity, that was mentioned in \cite{mk-mw_10} and can lead to rich families of MMOs is  \textit{Folded Saddle Node of type I} (FSN I). In the case of this bifurcation, small oscillations seem to be related to the presence of a delayed Hopf bifurcation rather than a true Hopf bifurcation.
For a more comprehensive overview we refer the reader to the recent review article \cite{md-jg-bk-ck-ho-mw_12}.
The theory of two slow and one fast variables has recently been generalized by Wechselberger \cite{mw_11} to arbitrary finite dimensions.  

In this article we study systems of two weakly coupled slow/fast oscillators. We assume that in the 
absence of the coupling one of the oscillators is undergoing a canard explosion and the other is
in general position. 
We show that turning on the weak coupling leads to the occurrence of
MMOs. We focus on the very interesting case of
FSN II, which, in the uncoupled case, corresponds to one of the oscillators undergoing a
canard explosion, while the other is at a stable equilibrium.  
As elaborated in \cite{mb-mk-mw_06} canard induced MMOs arise through a combined presence
of a folded singularity and suitable return mechanism. In this paper we focus on the existence
and properties of a folded singularity leaving the study of the return mechanism for future
investigations. 

The article is organized as follows. We start with a background section, Section \ref{sec-back}, in which we explain the 
very standard case of MMOs in the context of two slow and one fast variables and subsequently  the case
of two slow and two fast variables with a simple fold line. This material is included for completeness
and presented  in such a way that the coupled oscillator case becomes a simple
corollary.  In Section \ref{cos}, we treat the general case of coupled oscillators and 
in Section \ref{sec-ex}, we consider the example of two coupled FitzHugh-Nagumo systems.
\section{Background}\label{sec-back}
\subsection{The basic case of  one fast and two slow variables}\label{sec-back3D}
We consider a system of the form
 \begin{eqnarray}
    \eps\dot x &=&  f(x, y, \epsilon), \label{eq-gen3D_x}\label{back1}\\
    \dot y&=& g(x, y, \epsilon), \label{eq-gen3D_y}\qquad x\in \R,\; y=(y_1,y_2)\in\R^2, g=(g_1,g_2).\label{back2}
  \end{eqnarray}
The associated reduced system is 
\begin{eqnarray}
    0 &=&  f(x, y,0), \label{eq-gen3Dr_x}\\
    \dot y&=& g(x, y,0), \label{eq-gen3Dr_y}\qquad x\in \R,\; y=(y_1, y_2) \in\R^2.
  \end{eqnarray}
Let $S_0$ denote the surface defined by $0=f(x,y,0)$.
Non-hyperbolic points correspond to the points on $S_0$
for which $f_x=0$
(we use the notation $f_{\xi}=\partial f/\partial\xi$). At such points the equation $0=f(x,y,0)$ cannot be solved
for $x$ as a function of $y$.
Suppose $(0,0)$ is a non-hyperbolic point.
We make a non-degeneracy assumption $f_{xx}(0)\neq 0$.
In order to obtain an explicit equation for
the slow flow, we try to solve $0=f(x,y,0)$ for $y_1$ or $y_2$.
It is convenient to first change the variables to simplify the process
of finding such a solution.
We begin with a change of variables of the form $x \to x+\eta(y)$, $\eta(0,0)=0$,  where $\nu$ satisfies $f_x(\nu(y),y)=0$.
In the new variables, with additional scaling,  $f(x,y,0)$ has the form
\[
f(x,y_1,y_2)=f(0,y,0)+x^2+O( x^2|y|, x^3).
\]
In this last equation and the following, we do not write the variable $\epsilon$ in the function $f$, wich is always equal to $0$, i.e. we write  $f(x,y_1,y_2)$ for $f(x,y_1,y_2,0)$  .
We make a non-degeneracy assumption $f_y(0,0)\neq (0,0)$. Without lost of generality (WLOG), we can assume that $f_{y_1}(0,0)\neq 0$. We can now introduce
a new coordinate $\tilde y_1=-f(0,y)$, and immediately drop the tilde,
to simplify the notation.  In the new coordinates $f$ has the form
\[
f(x,y_1,y_2)=-y_1+x^2+O( x^2|y|, x^3).
\]
The transformations we made may be just
valid locally, i.e. only in a small neighborhood of the non-hyperbolic point.
Now $S_0$ is (locally) represented as a graph:
\begin{equation}\label{S03}
y_1=x^2(1+O(y_2, x)).
\end{equation}
and the fold is the straight line $x=y_1=0$.\\
We now transform \eref{eq-gen3Dr_x}-\eref{eq-gen3Dr_y}, first removing the constraint $f(x,y)=0$ and subsequently desingularising the resulting equation.
Differentiating \eref{S03} we get
\begin{equation}\label{S03p}
\dot y_1=2x(1+O(y_2, x))\dot x+O(x^2).
\end{equation} 
Substituting \eref{S03p} into \eref{eq-gen3Dr_y} we get
\begin{eqnarray}
2x(1+O(y_2, x))\dot x&=&g_1(x,x^2(1+O(y_2, x)),y_2)+O(x^2)\label{red3Dexp1}\\
\dot y_2 &=& g_2(x,x^2(1+O(y_2, x)),y_2).\label{red3Dexp2}
\end{eqnarray}
Now it is clear that \eref{red3Dexp1}  is singular at the fold
as long as $g_1(0,0,y_2)\neq 0$. Points on the fold line
for which $g_1(0,0,y_2)=0$
are called folded singularities. We would like to understand
the flow near such points better and to this end we 
apply a singular time rescaling to
\eref{red3Dexp1}-\eref{red3Dexp2} by multiplying the right-hand side (RHS) by the
factor $2x(1+O(y_2, x))$ and canceling it in \eref{red3Dexp1}.
This leads to the following system:
\begin{eqnarray}
\dot x&=&g_1(x,x^2(1+O(y_2, x)),y_2)+O(x^2)\label{reddes3D1}\\
\dot y_2 &=& 2x(1+O(y_2, x))g_2(x,x^2(1+O(y_2, x)),y_2).\label{reddes3D2}
\end{eqnarray}
We refer to \eref{reddes3D1}-\eref{reddes3D2} as the desingularized system.
Note that folded singularities correspond to equilibria of
\eref{reddes3D1}-\eref{reddes3D2} with $x=0$.
Further note that the trajectories of  
\eref{red3Dexp1}-\eref{red3Dexp2}
and \eref{reddes3D1}-\eref{reddes3D2}
restricted to the half plane $x>0$
differ only by time parametrization
but are the same as sets.
The trajectories in the half-plane $x<0$, on the attracting part of the critical manifold, 
are the same as sets but have the opposite
direction of time.
The flow near folded singularities is determined 
by the linearization of  \eref{reddes3D1}-\eref{reddes3D2}
at folded singularities, which is given by the Jacobian matrix:
\begin{equation}\label{jac3D}
\left (
\begin{array}{cc} 
g_{1,x}(0,0,y_2)&g_{1,y_2}(0,0,y_2)\\
2g_2(0,0,y_2)& 0
\end{array}
\right )
\end{equation}
Folded singularities are classified
according to the type of the corresponding equilibrium of 
\eref{reddes3D1}-\eref{reddes3D2}. Canards arise near folded saddles
and folded nodes and small oscillations are associated with folded
nodes, see the work of Beno\^it \cite{eb_90}, and 
Szmolyan and Wechselberger \cite{ps-mw_01}. 
Folded node, which is the only generic folded singularity
whose dynamics is accompanied by small oscillations, occurs  
when the corresponding equilibrium of  \eref{reddes3D1}-\eref{reddes3D2}
has two real positive eigenvalues (recall that we have changed the direction of the flow on the attracting side of the critical manifold).
\subsection{Folded singularities}\label{sec-foldsecs}
We now consider the case when $(f, g)$ (the RHS of \eref{back1}-\eref{back2}) depend on
a regular parameter and $g_2(0,0,y_2)$ passes through $0$ when this parameter is varied. Suppose that the equation $g_2(0,0,y_2)=0$ admits a unique solution $\bar{y}_2$.
If we assume that $g_{1,y_2}(0,0,y_2) neq 0$ for $y_2$ in a neighborhood of $\bar{y}_2$ and $\frac{\partial}{\partial y_2}g_2(0,0,\bar{y}_2)\neq 0$, this passage 
corresponds to a change of sign of the determinant
of the Jacobian \eref{jac3D}, from negative to positive or vice-versa. 
We make the asumption that $g_{1,x}>0$ wich guarantees that the flow of \eref{red3Dexp1}-\eref{red3Dexp2} on the stable part of the critical manifold is towards the fold.
This, in turn, means that the eigenvalues change from 
one positive and one negative (negative determinant) to two positive (positive determinant),
i.e. from saddle to node. It follows that this transition corresponds to the onset of
small oscillations. Hence the  onset of small oscillations is a consequence of the passage of $g_2(0,0,y_2)$ through $0$.
This transition is called Folded Saddle-Node of type II (FSN II).

As finding FSN II in coupled oscillator systems is the focus of this article we review some of
the features of this bifurcation, referring the reader to \cite{mk-mw_10} for details. An important feature of FSN II 
is that it corresponds to the passage of a true
equilibrium of \eref{eq-gen3Dr_x}-\eref{eq-gen3Dr_y} through the fold line. 
We assume WLOG that FSN II corresponds to the point
$(x,y,\lambda)=(0,0,0)$, where $\lambda$ is the regular parameter.
Suppose  that FSN II is {\em non-degenerate}, that is $(x,y,\lambda)=(0,0,0)$
corresponds to a transcritical bifurcation of  
\eref{reddes3D1}-\eref{reddes3D2}.
Then, for $\lambda=0$,
there is a folded singularity satisfying $g_2(0,0,0)=0$ and,
for $\lambda\neq 0$ but close to $0$ there is a point 
close to the origin, such that
$g_1(x,x^2+ O(x^3),y_2)=g_2(x,x^2+ O(x^3),y_2)=0$, i.e.
a true equilibrium of \eref{eq-gen3Dr_x}-\eref{eq-gen3Dr_y} near the fold line.
For the original system \eref{back1}-\eref{back2} one can prove that this corresponds
to a so called singular Hopf bifurcation \cite{jg_08}, and, in a different regime
of the parameter $\lambda$, to a delayed Hopf bifurcation \cite{mk-mw_10}.
These two bifurcations lead to the existence of small oscillations and thus enable
the existence of MMOs.
We note here that the interesting case of FSN II is when the equilibrium is stable on the
stable slow manifold and a saddle on the unstable slow manifold.

For folded nodes, not close to either FSN II
or to Folded Saddle Node of type I (FSN I), defined by $g_{1,y_2}=0$,
small oscillations are extremely small, which means that they
cannot be seen, except with very detailed numerics, see \cite{mb-mk-mw_06}. 
Hence the most interesting cases occur near FSN I or FSN II.
\subsection{Canonical system in two fast and two slow dimensions}\label{sec-back4D}
We consider a system with  two fast and two slow variables
\begin{eqnarray}
    \eps\dot x&=&  f(x, y,\eps), \label{eq-gensl_x}\\
    \dot y &=& g(x, y,\eps), \label{eq-gensl_y}\qquad x=(x_1,x_2)\in \R^2,\; y=(y_1,y_2)\in\R^2,
  \end{eqnarray}
  $f=(f_1,f_2)$, $g=(g_1,g_2)$.
 Recall that in the fast formulation \eref{eq-gensl_x}-\eref{eq-gensl_y} has the form
  \begin{eqnarray}
    x' &=&  f(x, y,\eps), \label{eq-gen_x}\\
    y' &=& \eps g(x, y,\eps) \label{eq-gen_y}.
  \end{eqnarray}
Recall also that the reduced problem has the form
 \begin{eqnarray}
    0 &=&  f(x, y,0), \label{eq-red_x}\\
    \dot{y} &=& g(x, y,0) \label{eq-red_y}.
  \end{eqnarray}
and the layer problem has the form
\begin{eqnarray}
    x' &=&  f(x, y,0), \label{eq-genlay_x}\\
    y' &=& 0\label{eq-genlay_y}.
  \end{eqnarray} 
Let $\Phi: \R^4\to\R^3$ be the map defined by
\begin{equation}\label{eq-defPhi}
\Phi(x,y)=\left (\begin{array}{c} f_1(x,y)\\f_2(x,y)\\\det(f_x)(x,y)
                       \end{array}\right ).
\end{equation}
The fold curve is defined by the condition $\Phi=0$.
We consider a point on the fold curve  (WLOG we assume that this point is the origin $(0,0)$).
We now state a few conditions which assure that the fold curve is {\it simple}; naturally
the first condition is that the linearization of layer system has a simple eigenvalue $0$. More specifically,
our first assumption is as follows:
\begin{eqnarray}\label{sfc1}
\hspace{-2 cm} \mbox{$f_x(0,0)$ has one simple $0$ eigenvalue and one simple eigenvalue in the left half-plane.}
\end{eqnarray}
Recall that in Section \ref{sec-back3D} we made additional assumptions,
namely $f_{xx}(0,0)\neq 0$ and $f_y(0,0)\neq (0,0)$. Here, we generalize these conditions
in the following way. Let $(v_1, v_2)$
denote the nullvector (the eigenvector of $0$) .
Our additional non-degeneracy conditions is as follows:
\begin{equation}\label{eq-nondeg4D}
D\Phi(0,0) \mbox{ has full rank and } \kappa=( \det(f_x))_x(0,0) \cdot (v_1, v_2)\neq 0.
\end{equation}
We have the following lemma:
\begin{Lemma}\label{lem-red3D}
Hypothesis \eref{eq-nondeg4D} implies that either $ \Phi_{(x_1, x_2, y_1)}(0,0)$ or $ \Phi_{(x_1, x_2, y_2)}(0,0)$
is invertible.
\end{Lemma}

 \noindent {\bf Proof}  Let $(w_1, w_2)$ the eigenvector of the negative eigenvalue.
 Note that either $f_{y_1}(0,0)$ or $f_{y_2}(0,0)$ must be linearly independent of $(w_1, w_2)$,
 otherwise $D\Phi(0,0)\R^4$ would be in the span of  $(w_1, w_2,0)$ and $(0,0,1)$.
 WLOG we assume that $f_{y_1}(0,0)$ is independent of $(w_1, w_2)$ and let 
 $M=\Phi_{(x_1, x_2, y_1)}(0,0)$.
 Note that $M(v_1, v_2, 0)=(0, 0, \kappa)$ and that the first two coordinates of $M(w_1, w_2, 0)$
are a multiple of $(w_1, w_2)$.  Note also that the first two coordinates of
 $M(0,0,1)$ are equal to $f_{y_1}(0,0)$. Since $f_{y_1}(0,0)$ 
 is independent of $(w_1, w_2)$, hence the vectors
 $M(v_1, v_2, 0)$,  $M(w_1, w_2, 0)$ and $M(0,0,1)$ are linearly independent. 
  
We now have the following proposition.
\begin{Proposition}\label{prop-sfc}
  We assume that $(0, 0)$ is on the fold curve, i. e.  $\Phi(0, 0)=0$, and that 
  conditions \eref{sfc1}-\eref{eq-nondeg4D}
  are satisfied. 
  Then, there exists a neighborhood $U=V\times W $ of $(0,0)\in \R^4$,  and a unique function $\psi:V\subset \R \rightarrow \R^3$, 
  such that: \[\{(x,y)\in U; \Phi=0\} = \{ (x,y_1)\in\R^3; (x,y_1)=\psi(y_2), y_2\in V\}.\]  
  \end{Proposition}
  
 \noindent {\bf Proof} By lemma \eref{lem-red3D} the matrix  $M=\Phi_{(x_1, x_2, y_1)}(0,0)$ is non singular.
Hence the equation $\Phi=0$ can be solved by the implicit function theorem, giving 
 $(x_1,x_2,y_1)$ as a function of $y_2$.  
 
In the remainder of this section we will assume that the fold curve has been transformed to 
 the coordinate line $(0,0,0,y_2)$ and that the fast system has been
diagonalized. More specifically, we assume that $f$ satisfies the following conditions
 \begin{eqnarray}
\hspace{-1.5 cm}f(0,0,0,y_2,0)=0,\quad \frac{\partial f_1}{\partial x_1}(0,0,0,y_2,0)=0, \quad \frac{\partial f_2}{\partial x_2}(0,0,0,y_2,0)=\lambda(y_2)<0, \label{eq-defcon1}\\
\frac{\partial f_1}{\partial x_2}(0,0,0,y_2,0)= \frac{\partial f_2}{\partial x_1}(0,0,0,y_2,0)=0. \label{eq-defcon2}
\end{eqnarray}
Note that the non-degenracy condition \eref{eq-nondeg4D} reduces to
\begin{equation}\label{eq-nondegcon1}
 \frac{\partial f_1}{\partial y_1}(0,0,0,y_2,0)\neq 0,\qquad \frac{\partial^2 f_1}{\partial x_1^2}(0,0,0,y_2,0)\neq 0.
 \end{equation}
WLOG we assume 
 \[
 \frac{\partial f_1}{\partial y_1}(0,0,0,y_2,0)=-1,\quad  \frac{\partial^2 f_1}{\partial x_1^2}(0,0,0,y_2,0)=2.
 \]
 We can now expand $f$ in Taylor series:
 \begin{equation}\label{eq-fexpts}
 f(x, y,\eps)=\left (\begin{array}{c} -y_1+x_1^2 + O(x_1x_2, x_2^2,x_1y_1,x_2y_1,y_1^2)+ O(||(x,y)||^3)+ O(\eps)\\
 \lambda(y)x_2+ O(y_1,\|x\|^2)+O(\eps).
 \end{array} 
 \right )
 \end{equation}
 
 The defining conditions of the slow manifold are $f(x,y,0)=0$, or
 \begin{eqnarray}
    0&=& -y_1+x_1^2 + O(x_1x_2, x_2^2,x_1y_1,x_2y_1,y_1^2)+ O(||(x,y)||^3), \label{eq-slow1} \\
    0 &=& \lambda(y)x_2+ O(y_1,\|x\|^2)\label{eq-slow2} 
  \end{eqnarray}
 From  \eref{eq-slow2} we get, by the implicit function theorem, $x_2=p(x_1,y)$, with $p(x_1,y)=O(y_1,x_1^2)$.
 Plugging into \eref{eq-slow1} we get the usual condition $y_1=x_1^2+ O(x_1^3)$.
 Following the approach of Section \ref{sec-back3D} we substitute $y_1=x_1^2+ O(x_1^3)$ into \eref{eq-red_y}
 getting
   \begin{eqnarray}
    (2x_1+O(x_1^2))\dot x_1 &=&  g_1(x_1, p(x_1,x_1^2+O(x_1^3),y_2),x_1^2+O(x_1^3),y_2,0), \label{eq-red1_x}\\
    \dot{y_2} &=& g_2(x_1, p(x_1,x_1^2+O(x_1^3),y_2),x_1^2+O(x_1^3),y_2,0). \label{eq-red1_y}
  \end{eqnarray}
  Finally we get the desingularized equation as follows:
  \begin{eqnarray}
    \dot x_1 &=&  g_1(x_1, p(x_1,x_1^2+O(x_1^3),y_2),x_1^2+O(x_1^3),y_2,0), \label{eq-des_x}\\
    \dot{y_2} &=& g_2(x_1, p(x_1,x_1^2+O(x_1^3),y_2),x_1^2+O(x_1^3),y_2,0)(2x_1+O(x_1^2)) \label{eq-des_y}
  \end{eqnarray}
  Folded singularities are points $(0, y_2)$ with $y_2$ satisfying $g_1(0,0,0,y_2,0)=0$, or, equivalently,
  equilibrium points of \eref{eq-des_x}-\eref{eq-des_y} on the fold line, i.e. satisfying $x_1=0$.
  The type of folded singularity is determined by the Jacobian of 
  \begin{equation}\label{jac4D}
\left (
\begin{array}{cc} 
g_{1,x_1}(0,0,0,y_2)&g_{1,y_2}(0,0,0,y_2)\\
2 g_2(0,0,0,y_2)& 0
\end{array}
\right ).
\end{equation}
Folded saddle nodes (FSN I and FSN II) can arise in the context of systems with two slow and two fast
dimensions in an analogous way as described in Section \ref{sec-foldsecs} for systems with two slow and one fast
dimensions, and, in the same
manner, correspond to the onset of MMOs.
\begin{Remark}
We want to emphasize here, link between the background presented in paragraph 2.1, which deals with 3d system and the computations done above, for the 4d system. As we have assumed in \eref{eq-defcon1} that $\lambda(y_2)<0$, we are able to write $x_2$ as a function of $O(x_1^2)$ by using the implicit function theorem in \eref{eq-slow1}-\eref{eq-slow2}. By this way, we can obtain a desingularized system in the 4d case, analog to the one obtained in the 3d case. This operation can be made in the same way, in a system with $N$ fast variables if we assume sufficiently hypothesis on eingenvalues of the jacobian matrix of  the fast subsystem on the fold line.
\end{Remark}
\section{Coupled oscillator system}\label{cos}

\subsection{Introduction of the system}
\label{cos1}
We consider a system of coupled oscillators in the following form
\begin{eqnarray}
\eps \dot x_{1}&=& -y_1+F_1(x_1)+\nu H^f_1(x,y)\label{eq-cosx1}\\
\eps \dot x_{2}&=& -y_2+F_2(x_2)+\nu H^f_2(x,y)\label{eq-cosx2}\\
\dot y_1 &=&G_1(x_1,y_1,c_1)+\nu H^s_1(x,y)\label{eq-cosy1}\\
\dot y_2 &=&G_2(x_2,y_2,c_2)+\nu H^s_2(x,y)\quad x=(x_1,x_2)\in \R^2,\; y=(y_1,y_2)\in \R^2\label{eq-cosy2}.
\end{eqnarray}
The parameters $\eps$ and $\nu$
are the singular parameter and the coupling parameter, respectively,
and are considered to be small. The parameters $c_1,c_2$
control the state of the uncoupled oscillators (moves the nullclines).
The coupling functions
$H^f:\R^4\to\R^2$, $H^s:\R^4\to\R^2$ and $F:\R^2\to\R^2$
are defined by $H^f=(H^f_1, H^f_2)$,
 by $H^s=(H^s_1, H^s_2)$ and $F=(F_1,F_2)$, respectively.
We assume that $y_j=F_j(x)$ are $S$ shaped curves.
Written on the fast time scale \eref{eq-cosx1}-\eref{eq-cosy2}  has the form:
\begin{eqnarray}
 x_{1}'&= -y_1+F_1(x_1)+\nu H^f_1(x,y)\label{eq-cofx1}\\
x_{2}'&= -y_2+F_2(x_2)+\nu H^f_2(x,y)\label{eq-cofx2}\\
y_1' &=\eps G_1(x_1,y_1,c_1)+\epsilon\nu H^s_1(x,y)\label{eq-cofy1}\\
y_2' &=\eps G_2(x_2,y_2,c_2)+\epsilon\nu H^s_2(x,y)\label{eq-cofy2}
\end{eqnarray}
We now find the conditions for the existence of a simple fold curve, as in
Section \ref{sec-back4D}. Let $\Phi$ defined as in section \ref{sec-back4D}.
 First note that the critical manifold of \eref{eq-cosx1}-\eref{eq-cosy2}, is defined by
 \begin{equation}\label{defcritman}
- y_j+F_j(x)+\nu H_j^f(x,y)=0,\qquad j=1,2.
 \end{equation}
 The linearization of the RHS of \eref{eq-cosx1}-\eref{eq-cosx2} is given by
 \begin{equation}\label{jacobgen}
\pmatrix{
 F'_1(x_1) +\nu H^f_{1,x_1}&\nu H^f_{1,x_2}\cr
\nu H^f_{2,x_1}&F_2'(x_2)+\nu H^f_{2,x_2}
}
\end{equation}
where we assume that $(x,y)$ is on the critical manifold,
i.e. satisfies \eref{defcritman}. The determinant
of  the matrix in \eref{jacobgen} 
is as follows:
\begin{equation}\label{foldet}
F'_1(x_1)F'_2(x_2)+\nu (F'_1(x_1)H^f_{2,x_2}+F'_2(x_2)H^f_{1,x_1})+ \nu^2\det H^f_x.
\end{equation}
%
 \begin{Proposition}\label{prop-sfc2}
Let $(x_0,y_0)$ satisfy $\Phi(x_0,y_0)=0$ for $\nu=0$. 
 We assume that
  \[F'_1(x_{0,1})=0 \mbox{ and } F'_2(x_{0,2})<0 \mbox{ as already specified in \eref{eq-defcon1}}.\] 
  And that,
  \[F_1''(x_{0,1})\neq 0.\]
  Then, there exists a neighborhood $U=V\times W $ of $(x_0,y_0)\in \R^4$,  and a unique function $\psi:V\subset \R \rightarrow \R^3$, such that: \[\{(x,y)\in U; \Phi=0\} = \{ (x,y_1)\in\R^3; (x,y_1)=\psi(y_2), y_2\in V\},\]  
   with $\psi(y_{2,0})=(x_{1,0},x_{2,0},y_{1,0})$. Furthermore, the parametrization of the fold curve
   depends also smoothly on the parameters $\nu$.  
  \end{Proposition}
   \noindent {\bf Proof}

We can apply Proposition \ref{prop-sfc} but we give here a direct proof by application of the  implicit function theorem. 
Let $M=\Phi_{(x_1, x_2, y_1)}|_{x=0,y=0,\nu=0}$.
Then, $\det M=F''(x_{10})(F'(x_{20}))^2\neq 0$. 
Note that, because the fast system does not depend on parameters $c_1,c_2$, the fold curve also does not.
 
 \subsection{Folded singularities and their nature}\label{sec-fsan}
 Let $x_1^*=\psi_1(y_2,\nu)$, $x_2^*=\psi_2(y_2,\nu)$
 and $y_1^*=\psi_3(y_2,\nu)$.
 We transform the fold curve to the line $(0,0,0,y_2)$, using the following change of variables:
 $\tilde{x_1}=x_1-x_1^*, \tilde{x_2}=x_2-x_2^*,\tilde{y_1}=y_1-y_1^*$. 
 We omit the tilde to simplify the notation. System \eref{eq-cosx1}- \eref{eq-cosy2}
 becomes
 \begin{eqnarray}
 \eps \dot x_{1}&=& -y_1+f_1(x_1,y_2)+\nu h^f_1(x,y)+O(\eps)\label{eq-costx1}\\
\eps \dot x_{2}&=& f_2(x_2,y_2)+\nu h^f_2(x,y)+O(\eps)\label{eq-costx2}\\
\dot y_1 &=&G_1(x_1+x_1^*,y_1+y_1^*,c_1)+\nu H^s_1(x+x^*,y_1+y_1^*, y_2)\nonumber\\
&&-\frac{d y_1^*}{d y_2}(G_2(x_2+x_2^*,y_2,c_2)+\nu H^s_2(x+x^*,y_1+y_1^*, y_2))\label{eq-costy1}\\
\dot y_2 &=&G_2(x_2+x_2^*,y_2,c_2)+\nu H^s_2(x+x^*,y_1+y_1^*, y_2),\label{eq-costy2}
\end{eqnarray}
where, 
\[f_1(x_1,y_2)=F_1'(x_1^*)x_1+F_1''(x_1^*)\frac{1}{2}x_1^2+O(x_1^3),\]
\[f_2(x_2,y_2)=F_2'(x_2^*)x_2+F_2''(x_2^*)\frac{1}{2}x_2^2+O(x_2^3),\]
and,
\[h^f_1(x,y)=H^f_{1,x_1}(x^*,y_1^*,y_2)x_1+H^f_{1,x_2}(x^*,y_1^*,y_2)x_2+H^f_{1,y_1}(x^*,y_1^*,y_2)y_1+O(||(x,y_1)||^2)\]
\[h^f_2(x,y)=H^f_{2,x_1}(x^*,y_1^*,y_2)x_1+H^f_{2,x_2}(x^*,y_1^*,y_2)x_2+H^f_{2,y_1}(x^*,y_1^*,y_2)y_1+O(||(x,y_1)||^2).\]
Our goal is  to arrive at the desingularized system \eref{eq-des_x}-\eref{eq-des_y}.
The first step is to diagonalize the linear part of the fast subsystem \eref{eq-costx1}-\eref{eq-costx2}. More precisely,
we change coordinates so that the RHS of \eref{eq-costx1}-\eref{eq-costy2} is transformed to the form \eref{eq-fexpts}.
Let 
\begin{equation}\label{egeivects}
\left (\begin{array}{c} 1\\v(y_2)\end{array}\right )\quad\mbox{and}\quad  
\left (\begin{array}{c} w(y_2)\\1 \end{array}\right )
\end{equation}
be the eigenvectors of the Jacobian \eref{jacobgen} at the points $(0,0,0,y_2)$
on the fold line, corresponding to the eigenvalues $0$ and
$\lambda(y_2)<0$, respectively. The diagonalizing transformation has the form
\[
\tilde x = P(y)^{-1}x,
\]
where the columns of the matrix $P$ are the eigenvectors \eref{egeivects}.
In the new variables system \eref{eq-costx1}-\eref{eq-costy2} becomes (tilde is omitted)

 \begin{eqnarray}
\eps\dot x_{1}=& \frac{1}{1-vw} \big(-y_1+ \frac{K}{2}x_1^2 +\nu(  O(y_1))+O(x_1x_2, x_2^2) +O(\eps) \big) \label{eq-cosdiax1}\\
\eps\dot x_{2}=& \lambda(y_2)x_2+\nu(  O(y_1))+O(||x||^2)+ O(\eps) \label{eq-cosdiax2}\\
\dot y_1 =& G_1(x_1+x_2w(y_2)+x_1^*,y_1+y_1^*,c_1)+\nu H^s_1(P(y_2)x+x^*,y_1+y_1^*, y_2)\nonumber\\
& \hspace{-1.5cm}  -\frac{d y_1^*}{d y_2}\Big (G_2(x_2+x_1v(y_2)+x_2^*,y_2,c_2) +\nu H^s_2(P(y_2)x+x^*,y_1+y_1^*, y_2)\Big ) 
&\label{eq-cosdiay1}\\
\dot y_2 =& G_2(x_2+x_1v(y_2)+x_2^*,y_2,c_2)+\nu H^s_2(P(y_2)x+x^*,y_1+y_1^*, y_2),\label{eq-cosdiay2}
\end{eqnarray}
where,
\begin{equation}
K=F_1''(x_1^*)+O(\nu).
\label{eq:K}
\end{equation}
WLOG, in the remaining of this section, we will assume that $F_1''(x_1^*)>0.$
Note that the fast subsystem \eref{eq-cosdiax1}-\eref{eq-cosdiax2} is now as specified in \eref{eq-fexpts}. 
Applying the procedure described in Section \ref{sec-back4D} we obtain the reduced system
\begin{eqnarray}
    (Kx_1+O(x_1^2))\dot x_1 &=  g_1(x_1, y_2,c_1,c_2), \label{eq-redcosc_x}\\
    \dot{y_2}&= g_2(x_1, y_2,c_2), \label{eq-redcosc_y}
  \end{eqnarray}
  where,
  \begin{eqnarray} 
 g_1(x_1,y_2,c_1,c_2)=&G_1(x_1^*+x_1+w(y_2)O(x_1^2),y_1^*+O(x_1^2),c_1)\nonumber\\
 & \hspace{-2cm}+\nu H^s_1(x_1^*+x_1+w(y_2)O(x_1^2), x_2^*+v(y_2)x_1+O(x_1^2),y_1^*+O(x_1^2), y_2)\nonumber\\
& \hspace{-2cm} -\frac{d y_1^*}{d y_2}\Big( G_2(x_2^*+v(y_2)x_1+O(x_1^2)) +\nu H^s_2(x_1^*+x_1+w(y_2)O(x_1^2),\nonumber\\
& x_2^*+v(y_2)x_1+O(x_1^2),y_1^*+O(x_1^2), y_2)\Big ) \label{eq-defg1}\\
g_2(x_1, y_2,c_2)=& G_2(x_2^*+v(y_2)x_1+O(x_1^2)) +\nu H^s_2(x_1^*+x_1+w(y_2)O(x_1^2),\nonumber\\ &x_2^*+v(y_2)x_1+O(x_1^2),y_1^*+O(x_1^2), y_2),\label{eq-defg2}
  \end{eqnarray}
and the desingularized system
\begin{eqnarray}
    \dot x_1 &=& g_1(x_1, y_2,c_1,c_2), \label{eq-descosc_x}\\
    \dot{y_2} &=&g_2(x_1, y_2,c_2)(2Kx_1+O(x_1^2)). \label{eq-descosc_y}
  \end{eqnarray}
\subsection{Folded singularities of type FSN II and the existence of MMOs}
\label{sec-FSN IImmo}
As discussed in Section \ref{sec-foldsecs} a well known mechanism 
of transition to MMOs is Folded Saddle-Node of type II (FSN II),
see \cite{mk-mw_10}, \cite{ScholarWech:2007} . This is a codimension one transition,
corresponding to the passage of the system from a parameter region
with an excitable equilibrium to a parameter region of MMO dynamics.
It can be described as follows in the context of system \eref{eq-redcosc_x}-\eref{eq-redcosc_y}: as the regular parameter is varied
a stable equilibrium of \eref{eq-redcosc_x}-\eref{eq-redcosc_y} approaches the fold, and, for the critical
value of the regular parameter, satisfies the conditions:
\begin{eqnarray}
 g_1(0, y_2,c_1,c_2)&=0 \label{eq-conFSN II1}\\
g_2(0,y_2,c_2)&=0\label{eq-conFSN II2}.
\end{eqnarray} 
On the other side of criticality there are MMOs as well as an equilibrium, now unstable.

\begin{figure}
\subfigure[]{\includegraphics[scale=0.3, angle=-90]{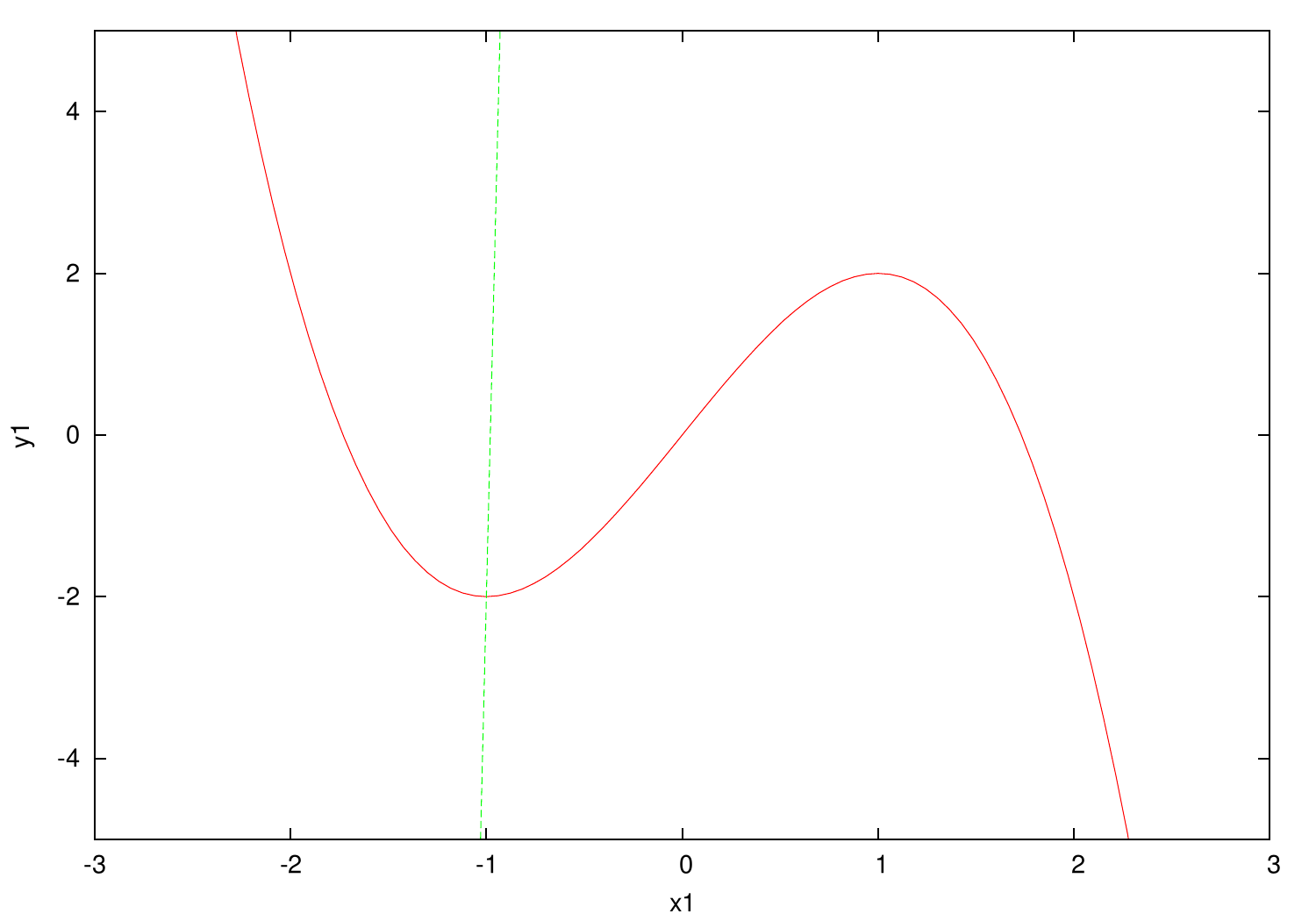}}
\subfigure[]{\includegraphics[scale=0.3, angle=-90]{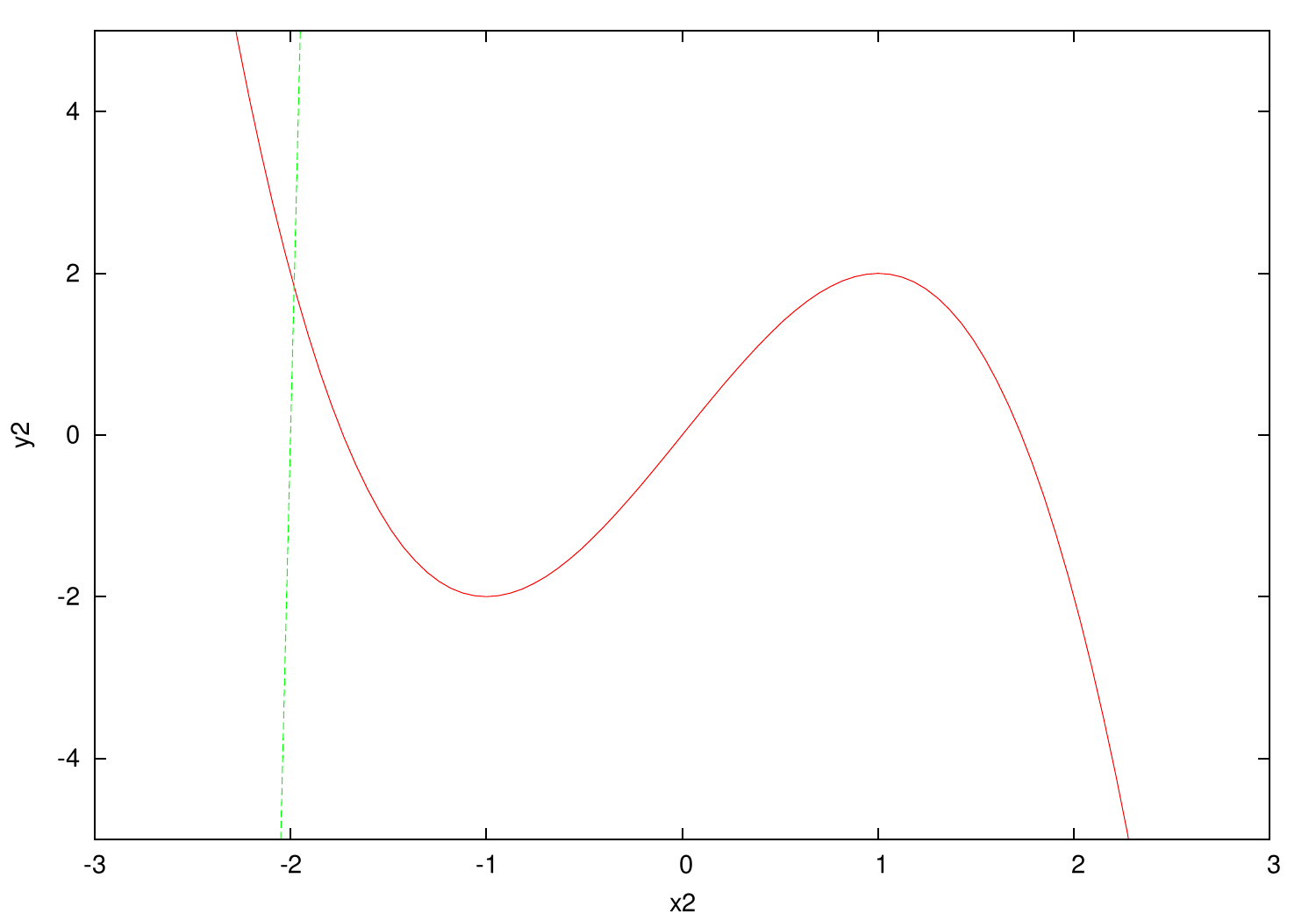}}
\caption{The $\nu=0$ limit configurations. Panel a) The uncoupled system \eref{eq-cofx1}-\eref{eq-cofy1} admits  a unique repulsive stationary point near the fold. Panel b) The uncoupled system \eref{eq-cofx2}-\eref{eq-cofy2} admits  a unique attractive stationary point on the stable manifold. }
\label{fig-se}
\end{figure}

In the following, we assume that the parameter $c_2$ is fixed. 
Note that 
for $\nu=0$, $y_1^*$ does not depend on $y_2$ and thus $dy_1^*/d y_2=O(\nu)$. We assume  that, for $\nu=0$, the equation $G_2(x_2^*,y_2,c_2)=0$ has a unique solution $x^*_{2}=\psi_2(y_{2},c_2)$  which is
a stable equilibrium of the uncoupled  $(x_2,y_2)$ subsystem, and that $G_{2,y_2}\neq 0$ at this point. By the implicit function theorem this gives a unique solution $y_2$ for the equation \eref{eq-conFSN II2}, for $\nu$ small enough. Recall that, by hypothesis, in the case $\nu=0$ we assume that, each uncoupled subsystem \eref{eq-cofx1}-\eref{eq-cofy1} and \eref{eq-cofx2}-\eref{eq-cofy2} admits a unique stationary point, one attractive on the stable manifold and one repulsive near the fold. This corresponds to the nullcline configuration shown in Figure \ref{fig-se}. 
Let us denote this particular value $\bar{y}_{2}$. It follows that this value  $\bar{y}_{2}$ determine particular values of $(x_1^*,y_1^*,x_2^*)=(\bar{x}_1^*,\bar{y}_1^*,\bar{x}_2^*)$. We assume that solving equation \eref{eq-conFSN II1} with these values gives a unique value of $c_1=\bar{c_1}$. 

Recall from the discussion in Section \ref{sec-foldsecs} that a non-degenerate FSN II singularity corresponds to 
a transcritical bifurcation for system \eref{eq-descosc_x}-\eref{eq-descosc_y}.

The following proposition establishes the existence of  a transcritical bifurcation for \eref{eq-descosc_x}-\eref{eq-descosc_y} and the existence of FSN II.
\begin{Proposition}
\label{PropFSN IIosc}
Let us assume that, for system \eref{eq-descosc_x}-\eref{eq-descosc_y}:
\begin{eqnarray}
\frac{dg_1}{dy_2}\neq 0 \mbox{ and }  \frac{dG_1}{dx_1} \frac{dG_2}{dy_2}\neq 0. \label{H1O}
\end{eqnarray}
Then for $c_1$ in a neighborhood of $c_1=\bar{c_1}$ and $\nu$ small enough, there exists two stationary points for system \eref{eq-descosc_x}-\eref{eq-descosc_y}:
$(x_{1,e},y_{2,e})$ and $(0,y_{2,fold})$.
Furthermore, if we assume that
\begin{eqnarray}
\frac{dg_1}{dy_2}\frac{dg_2(0,y_{2,fold}(c_1))}{dc_1}<0,
, G_{1,x_1}>0 \mbox{ and }   G_{2,y_2}\frac{dx_{1,e}}{dc_1}<0\label{H2O},
\end{eqnarray}
then, for  $\nu$ sufficiently small, there  is a transcritical bifurcation. As the parameter $c_1$ increases from left to right, the folded stationary point passes from  a  folded saddle to a repulsive folded node, whereas the ordinary stationary point passes from a repulsive node  to a  saddle. It follows that $(0,\bar{y}_2)$, is an FSN II point.  
\end{Proposition}
 \textbf{Proof}\\
The folded stationary point, $(0,y_{2,fold})$ is obtained by solving the equation $g_1(0,y_2,c_1,c_2)=0, i.e, g_1=0,x_1=0.$
For this, we apply the implicit function theorem to $g_1(0,y_2,c_1,c_2)$ as a function of $y_2$ and $c_1$. For
$y_2=\bar{y}_2$ and $c_1=\bar{c}_1$, we have $g_1=0$. From hypothesis \eref{H1O}, it follows that $\frac{dg_1}{dy_2}(0,\bar{y}_2,\bar{c}_1,c_2)\neq 0$, then the implicit function theorem gives the existence of a stationary point $(x_1=0,y_2=y_{2,fold})$ for $c_1$ in a neighborhood of $\bar{c}_1$. 
The second stationary point, $(x_{1,e},y_{2,e})$, is obtained by solving the equation $g_1(x_1,y_2,c_1)=0, g_2(x_1,y_2)=0.$ 
Note that this second stationary point corresponds to a true stationary point of system \eref{eq-cosx1}-\eref{eq-cosy2}. For $\nu=0$, we know that these equations admit a unique solution.  Using hypothesis \eref{H2O}, we can apply the implicit function theorem in a neighborhood of $\nu=0$.
This gives the existence of a unique stationary point for $\nu$ sufficiently small. 
For the first stationary point, the folded one $(0,y_{2,fold})$ the jacobian is given by
\begin{equation}\label{jacosc1}
\hspace*{-2cm}J_1=\left (
\begin{array}{cc} 
G_{1,x_1}+O(\nu)&\frac{dg_1}{dy_2}\\
2 K g_2& 0
\end{array}
\right ).
\end{equation}
This gives:
\[\det J_1=-(\frac{dg_1}{dy_2})2 K g_2.\]
For the second stationary point, the ordinary one, the jacobian is given by:
\begin{equation}\label{jacosc2}
J_2=\left (
\begin{array}{cc} 
G_{1,x_1}+O(\nu)+O(x_1)&\frac{dg_1}{dy_2}\\
KG_{2,x_1}x_1+O(x_1^2)& KG_{2,y_2}x_1+O(x_1^2)
\end{array}
\right ).
\end{equation}
This gives:
\[\det J_2=KG_{1,x_1}G_{2,y_2}x_1+O(x_1\nu).\]
From hypothesis \eref{H2O}, it follows that, as $c_1$ crosses the value $\bar{c}_1$, $\det  J_1$ crosses the value $0$ from negative to positive whereas $ \det J_2$ crosses the value $0$ from positive to negative. Subsequently the stationary point $(0,y_{2,fold})$ bifurcates from a saddle to a repulsive node whereas the stationary point $(x_{1,e},y_{2,e})$ bifurcates from a repulsive node to a saddle. Then, it follows that the system \eref{eq-descosc_x}-\eref{eq-descosc_y} admits a transcritical bifurcation at point  $(0,\bar{y}_2)$ when $c_1$ crosses the value $\bar{c}_1$, and $(0,\bar{y}_2)$ is a Folded Saddle-Node of type II.

\section{Example -- coupled FitzHugh-Nagumo systems.}
\label{sec-ex}

\subsection{Simple fold line}
We consider the following system
\begin{center}
\begin{equation}
\label{FHNdl}
 \left \{ 
     \begin{array}{rcl}
      \epsilon \dot{x}_1&=& F(x_1)-y_1+\alpha_1(x_2-x_1) \\
\epsilon \dot{x}_2&=& F(x_2)-y_2+\alpha_2(x_1-x_2) \\
      \dot{y}_1&=& x_1-c_1+\beta_1(y_2-y_1)\\
      \dot{y}_2&=& x_2-c_2+\beta_2(y_1-y_2)\\
     \end{array}
      \right. 
\end{equation}
\end{center}
with
\[F(z)=-z^3+3z.\]
We define
\begin{center}
\begin{equation}
\label{f}
 l(x_1,y_1,x_2,y_2)=\left \{ 
     \begin{array}{l}
      F(x_1)-y_1+\alpha_1(x_2-x_1) \\
      F(x_2)-y_2+\alpha_2(x_1-x_2) \\
     \end{array}
      \right. 
\end{equation}
\end{center}

Points on the fold line are defined by the condition
\begin{center}
\begin{equation}
\label{fold}
 \left \{ 
     \begin{array}{rcl}
    \det(Dl(x,y))&=& 0 \\
      l_1(x,y)&=& 0\\
l_2(x,y)&=& 0 \\
     \end{array}
      \right. 
\end{equation}
\end{center}
wich, for system \eref{FHNdl}, is given by
\begin{center}
\begin{equation*}
\label{systemfold}
 \left \{ 
     \begin{array}{rcl}
      (F'(x_1)-\alpha_1)(F'(x_2)-\alpha_2)-\alpha_1\alpha_2&=& 0 \\
      y_1&=& F(x_1)+\alpha_1(x_2-x_1)\\
y_2&=&  F(x_2)+\alpha_2(x_1-x_2).\\
     \end{array}
      \right. 
\end{equation*}
\end{center}

Thanks to proposition \ref{prop-sfc2}, the existence of a \textit{simple} fold line follows from $F'(x_1)=0$, $F'(x_2)\neq 0$ and $F''(x_1)\neq 0$
(the alternative choice would be $F'(x_2)=0$, $F'(x_1)\neq 0$ and $F''(x_2)\neq 0$).  
We  obtain a solution of system \eref{systemfold} with $(x_1,x_2,y_1)$ given as function of $y_2$.
We denote this solution by $(x_1^*,x_2^*,y_1^*)=(\varphi_1(y_2),\varphi_2(y_2),\varphi_3(y_2))$.
Note that we have obtained a curve in $\R^4$  parametrized by $y_2$, see \ref{prop-sfc2} for a general statement.

\subsection{Folded singularities and their nature}

System \eref{eq-costx1}-\eref{eq-costy2} in context of \eref{FHNdl} becomes
\begin{center}
\begin{equation}
\label{a}
 \left \{ 
     \begin{array}{rcl}
      \epsilon \dot{x}_1&=&-y_1+f_1(x_1,y_2)+\alpha_1(x_2-x_1)+O(\epsilon)\\
\epsilon \dot{x}_2&=&f_2(x_2,y_2)+\alpha_2(x_1-x_2)+O(\epsilon) \\
    \dot{y}_1&=& x_1^*+x_1-c_1+\beta_1(y_2-(y_1^*+y_1))\\
      & &-\frac{d y_1^*}{d y_2}(x_2^*+x_2-c_2+\beta_2(y_1^*+y_1-y_2))\\
     \dot{y}_2&=& x_2^*+x_2-c_2+\beta_2(y_1^*+y_1-y_2)\\
     \end{array}
      \right. 
\end{equation}
\end{center}
where
\[f_1(x_1)=F'(x_1^*)x_1+ F''(x_1^*)\frac{x_1^2}{2}-x^3_1\]
and
\[f_2(x_2)= F'(x_2^*)x_2+F''(x_2^*)\frac{x_2^2}{2}-x^3_2.\]
As in section \ref{sec-fsan}, we diagonalize the linear part of the fast system.
The eigenvalues of $Df(x_1^*,x_2^*$) are: 
\[0 \mbox{ and } \lambda(y_2)=F'(x_1^*)+F'(x_2^*)-\alpha_1-\alpha_2.\]
As explained in section \eref{sec-FSN IImmo}, we focus on dynamics decribed  by configuration in Figure \ref{fig-se}, thus $F'(x_2^*)<0$.
The associated eigenvectors are, respectively,\\
\vspace{0.5cm}
$\displaystyle P_1=\pmatrix{
 1\cr
\frac{\alpha_2}{\alpha_2-F'(x_2^*)}}
$ and $ \displaystyle P_2=
\pmatrix{
\frac{\alpha_1}{F'(x_2^*)-\alpha_2}\cr
 1
}.$\\
\vspace{0.5cm}
We now diagonalize the fast subsystem using the transformation
\[X=P\tilde{X} \mbox{ with } X=\pmatrix{
        x_1\cr
        x_2}
        \]  
where  $P$ is the matrix whose columns are the eigenvectors.
We introduce the following notations:
\[w=\frac{\alpha_1}{F'(x_2^*)-\alpha_2}, \quad v=\frac{\alpha_2}{\alpha_2-F'(x_2^*)} \quad d=\frac{1}{1-vw}.\] 
After some transformations , system \eref{a} becomes:
\begin{center}
\begin{equation*}
\label{b}
 \left \{ 
     \begin{array}{rcl}
      \epsilon \dot{x}_1&=& d((F''(x_1^*)\frac{(x_1+wx_2)^2}{2}-(x_1+wx_2)^3-y_1)\\
      & &-w (F''(x_2^*)\frac{(vx_1+x_2)^2}{2}-(vx_1+x_2)^3))+O(\epsilon)\\
\epsilon \dot{x}_2&=& \lambda(y_2)x_2+d(-v(F''(x_1^*)\frac{(x_1+wx_2)^2}{2}-(x_1+wx_2)^3-y_1)\\
& &+F''(x_2^*)\frac{(vx_1+x_2)^2}{2}-(vx_1+x_2)^3) +O(\epsilon)\\
 \dot{y}_1&=& x_1^*+x_1+wx_2-c_1+\beta_1(y_2-(y_1^*+y_1))\\
      & &-\frac{d y_1^*}{d y_2}(x_2^*+vx_1+x_2-c_2+\beta_2(y_1^*+y_1-y_2))\\
     \dot{y}_2&=& x_2^*+vx_1+x_2-c_2+\beta_2(y_1^*+y_1-y_2).\\
     \end{array}
      \right. 
\end{equation*}
\end{center}
Hence, the reduced equation is given by
\begin{eqnarray}
      0&=& d((F''(x_1^*)\frac{(x_1+wx_2)^2}{2}-(x_1+wx_2)^3-y_1)\\\nonumber
      &&-w (F''(x_2^*)\frac{(vx_1+x_2)^2}{2}-(vx_1+x_2)^3))\label{Reqx1}\\
0&=& \lambda(y_2)x_2+d(-v(F''(x_1^*)\frac{(x_1+wx_2)^2}{2}-(x_1+wx_2)^3-y_1)\\\nonumber
&&+F''(x_2^*)\frac{(vx_1+x_2)^2}{2} -(vx_1+x_2)^3) \label{Reqx2}\\
      \dot{y}_1&=& x_1^*+x_1+wx_2-c_1+\beta_1(y_2-(y_1^*+y_1)) \label{Reqy1}\\\nonumber
      & &-\frac{d y_1^*}{d y_2}(x_2^*+x_2-c_2+\beta_2(y_1^*+y_1-y_2)) \label{Reqy2}\\
     \dot{y}_2&=& x_2^*+vx_1+x_2-c_2+\beta_2(y_1^*+y_1-y_2)\\\nonumber
\end{eqnarray}
Now, using equation  \eref{Reqx1} and the implicit function theorem in \eref{Reqx2}, as $\lambda(y_2)<0$, one can obtain $x_2$ as a function of $x_1,y_2$
wich leads to, 
\begin{equation*}
 x_2= O(x_1^2).
\end{equation*}
Now, \eref{Reqx1} can be rewritten in the form
\begin{equation}\label{solReqx1}
 y_1= \frac{K}{2}x_1^2+O(x_1^3)
\end{equation}
with
\begin{eqnarray*}
K&=F''(x_1^*)-wv^2F''(x_2^*)\\
&=F''(x_1^*)+O(\nu^3),
\end{eqnarray*}
where $0\leq\nu\leq \max(\alpha_1,\beta_1,\alpha_2,\beta_2)$.
Subsequently, we derivate \eref{solReqx1}, and plug into \eref{Reqy1} and \eref{Reqy2}. We obtain,
\begin{eqnarray*}
      (Kx_1+O(x_1^2))\dot{x}_1&=& x_1^*+x_1-c_1+\beta_1(y_2-y_1^*)\\
      & &-\frac{d y_1^*}{d y_2}(x_2^*+bx_1-c_2+\beta_2(y_1^*-y_2))+O(x_1^2)\\
     \dot{y}_2&=& x_2^*+bx_1-c_2+\beta_2(y_1^*-y_2)+O(x_1 ^2)\\
\end{eqnarray*}
This gives the following desingularized system:
\begin{eqnarray}
     \dot{x}_1&=& x_1^*+x_1-c_1+\beta_1(y_2-y_1^*)\nonumber\\
     & &-\frac{d y_1^*}{d y_2}(x_2^*+vx_1-c_2+\beta_2(y_1^*-y_2))+O(x_1^2)\label{Deseq1} \\
       \dot{y}_2&=& (x_2^*+vx_1-c_2+\beta_2(y_1^*-y_2)+O(x_1 ^2))(Kx_1+O(x_1^2))\label{Deseq2}  
\end{eqnarray}
Hence using notations of section \ref{sec-FSN IImmo}, we have:
\[g_1(x_1,y_2)=x_1^*+x_1-c_1+\beta_1(y_2-y_1^*)-\frac{\partial y_1^*}{\partial y_2}(x_2^*+vx_1-c_2+\beta_2(y_1^*-y_2))+O(x_1^2)\]
and
\[g_2(x_1,y_2)=x_2^*+vx_1-c_2+\beta_2(y_1^*-y_2)+O(x_1 ^2).\]
\begin{Proposition}
\label{PropFNFHN}
Let us assume $c_2<-1$ and $
    \beta_1-\frac{\alpha_1}{9(1-c_2^2)^2}>0$. For $\alpha_1, \alpha_2,  \beta_1, \beta_2$ small enough,  let $\bar{y}_2$ be the solution of  $g_2(0,y_2)=0$ and define  $\bar{c}_1=x_1^*+\beta_1(\bar{y}_2-y_1^*)$.
Then, for $\alpha_1, \alpha_2, \beta_1, \beta_2$ in  a small neighborhood of zero,  the system \eref{Deseq1}-\eref{Deseq2}, admits a transcritical bifurcation at point $x_1=0,y_2=\bar{y}_2, c_1=\bar{c}_1$.  As the parameter $c_1$ crosses the value $\bar{c}_1$ from left to right, the folded stationary point passes from  a  folded saddle to a repulsive folded node, whereas the ordinary stationary point passes from a repulsive node  to a  saddle. 
\end{Proposition}
\noindent \textbf{Proof}\\
This proposition is a specific case of proposition \eref{PropFSN IIosc}. We will verify that:
\begin{eqnarray}
\frac{dg_1}{dy_2}(0,\bar{y}_2,\bar{c}_1)>0 \label{H1},\\
\frac{dG_1}{dx_1}>0, \label{H2}\\
\frac{dG_2}{dy_2} < 0,  \label{H3}\\
\frac{dx_{1,e}}{dc_1}>0, \label{H4}\\
 \frac{dg_2(0,\bar{y}_{2,fold}(c_1))}{dc_1}<0   \label{H5}.
\end{eqnarray}
We start with hypothesis \eref{H1}. After some computation, we find that this  hypothesis reads as,
\begin{eqnarray*}
    \beta_1-\frac{\alpha_1}{9(1-\bar{x}_2^{*2})^2}+O(\nu^2) &>0.\\   
\end{eqnarray*}
This holds since $ \bar{x}_2^{*}=c_2+O(\nu)$.
 Hypothesis \eref{H2}, is verified since $\frac{\partial G_1}{\partial x_1}=1$ .
 
 Now, we deal with \eref{H3}. Here, it reads: 
\begin{eqnarray}
\frac{dx_2^*}{dy_2} <0 . 
\end{eqnarray}
This holds since $\frac{dx_2^*}{dy_2}=\frac{1}{F'(x_2)}+O(\nu)<0$ for $\nu$ sufficiently small.
Hypothesis \eref{H4} is verified since we have $\frac{dx_{1,e}}{dc_1}=1+ O(\nu)>0$ .

Finally, we come to hypothesis \eref{H5}. From \eref{H1}, we know that $\frac{dg_1}{dy_2}(0,y_2)>0$. It follows that $\frac{y_{2,fold}}{dc_1}>0$. Furthermore, we have that  $\frac{dg_2}{dy_2}(0,y_2)<0$. We conclude that 
 $\frac{dg_2(0,\bar{y}_{2,fold}(c_1))}{dc_1}<0$.

\subsection{Numerical simulations}
We have performed numerical integration of system \eref{FHNdl} on time interval $[0, 500]$, using a fourth-order Runge-Kutta  method with time step of $0.0001$. The parameter values are the following: $\eps=0.01, \alpha_1=\beta_1=\alpha_2=\beta_2=0.05, c_2=-1.5$. As explained in proposition \ref{PropFNFHN}, we obtain the critical value of parameter $\bar{c}_1$  by solving the following equations:
\begin{eqnarray*}
g_2(0,y_2)&=0\\
c_1&=x_1^*+\beta_1(y_2-y_1^*).
\end{eqnarray*}
In our case, we find:
\[\bar{c}_1\simeq -0.95266.\]

\begin{figure}[ht!]
\begin{center}
\vspace{3cm}
\subfigure[]{\includegraphics[scale=.3,angle=-90]{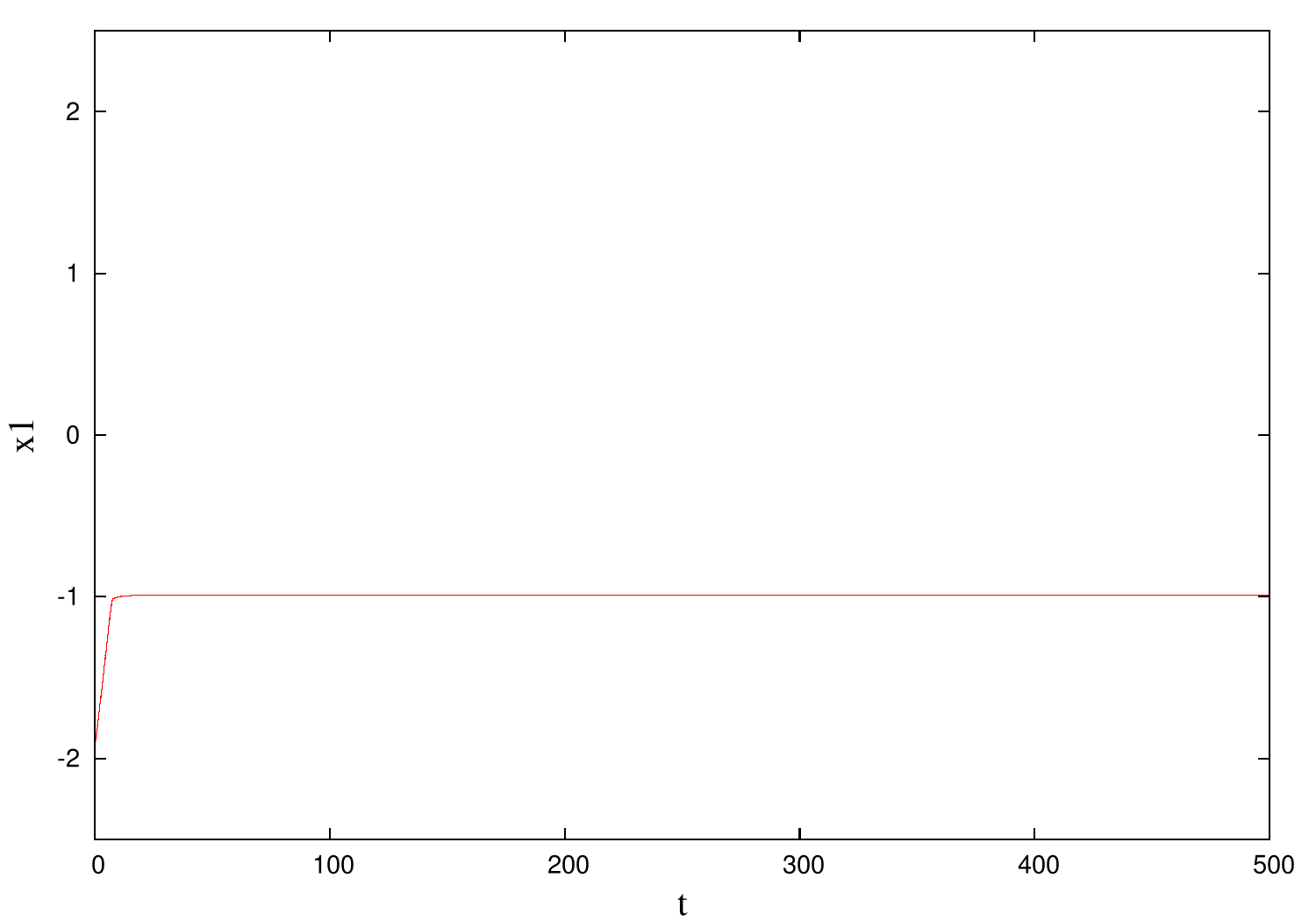}}
\subfigure[]{\includegraphics[scale=.3,angle=-90]{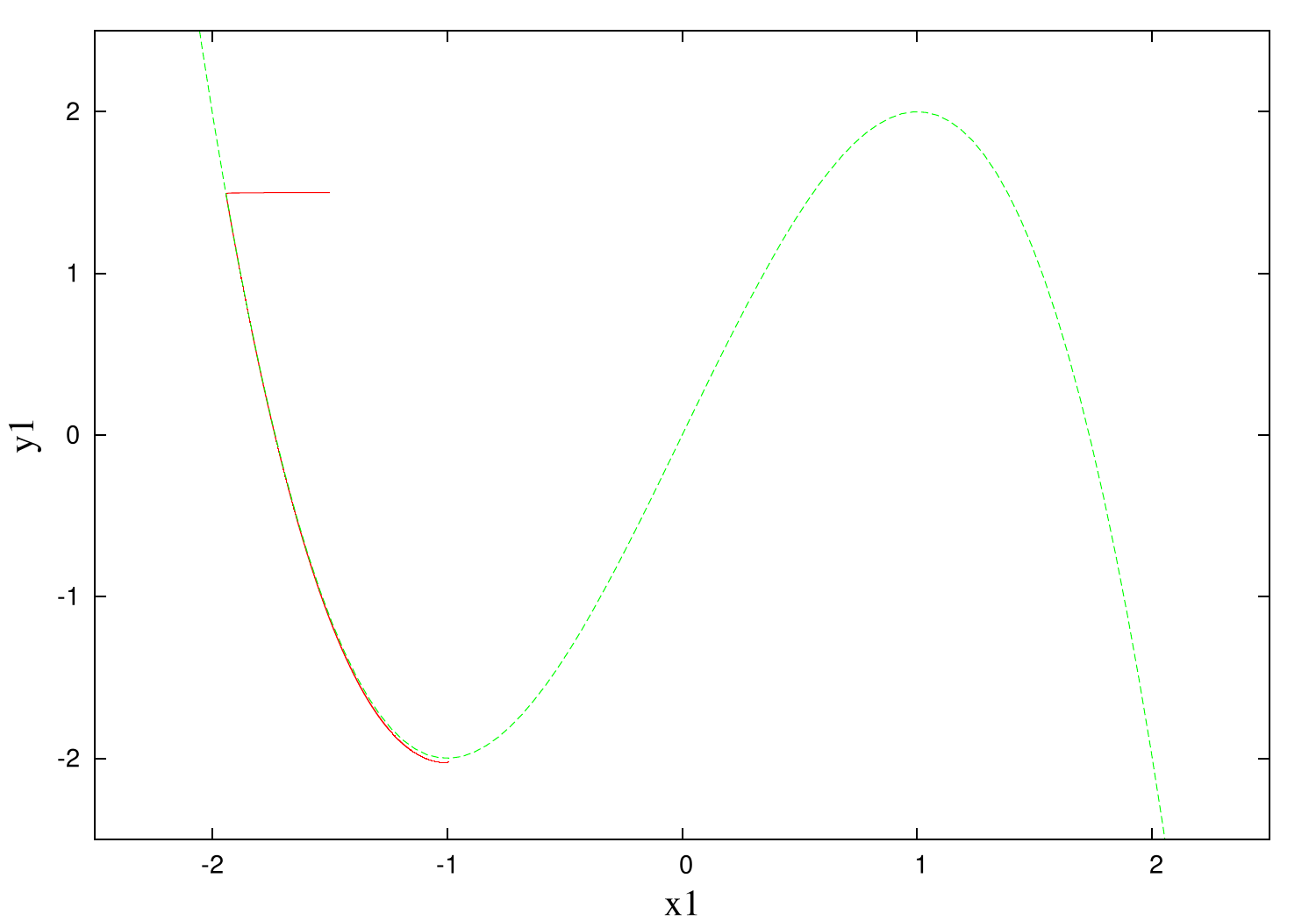}}\\
\caption{This figure corresponds to the simulation of system \eref{FHNdl} for  a value of parameter $c_1=-0.953$. In panel (a), we show the time-evolution of $x_1$ in the interval [0,500]. In panel (b), we show the behavior of the system  in the $(x_1, y_1)$ phase-plane. For this value of the parameter $c_1$ the system evolves asymptotically to his stationary state.}
\label{pointstable2Dfig}
\end{center}
\end{figure}

\begin{figure}[ht!]
\begin{center}
\subfigure[]{\includegraphics[scale=.3,angle=-90]{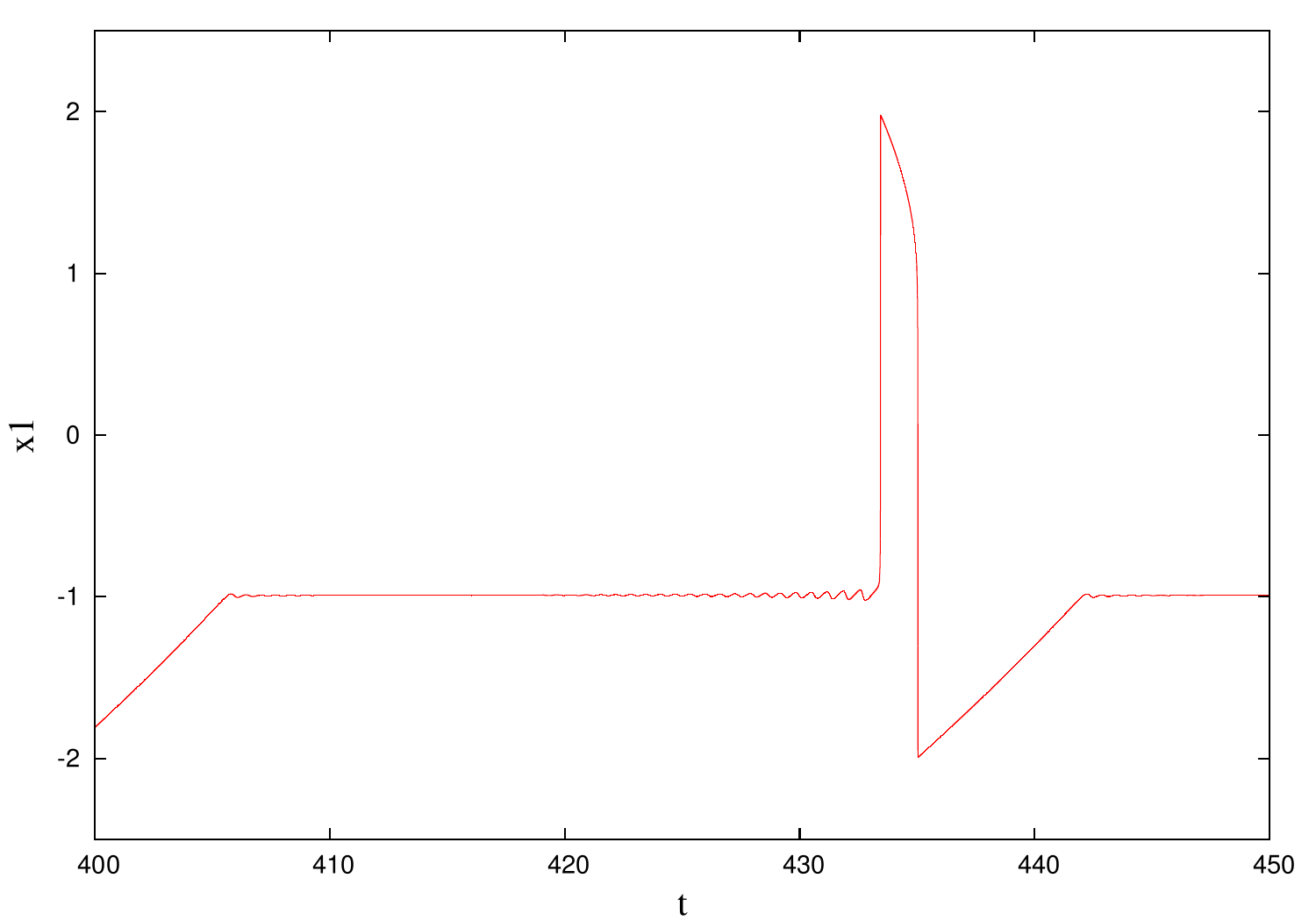}}
\subfigure[]{\includegraphics[scale=.3,angle=-90]{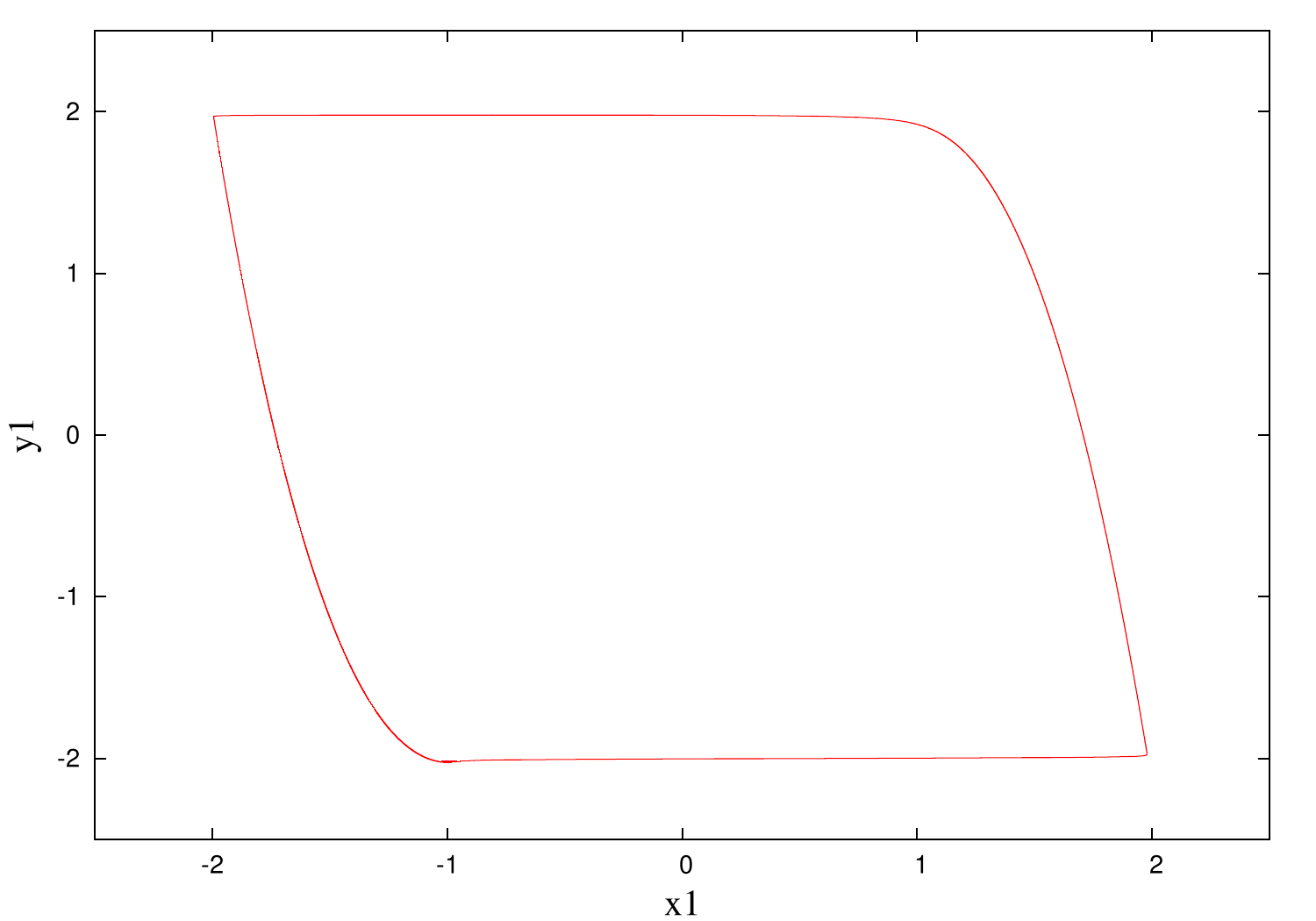}}\\
\caption{This figure corresponds to the simulation of system \eref{FHNdl} for a value of parameter $c_1=-0.952$. In panel (a), we show the time-evolution of $x_1$, for $t$ in the interval [400,450]. In panel (b), we show the behavior of the system  in the $(x_1, y_1)$ phase-plane. For this value of the parameter $c_1$, the MMOs appear. But the amplitude of the small oscillations is very small compared to that of large oscillations, and a zoom in necessary  to distinguish them.}
\label{plc2Dfig}
\end{center}
\end{figure}

\begin{figure}
\begin{center}
\subfigure[]{\includegraphics[scale=.3,angle=-90]{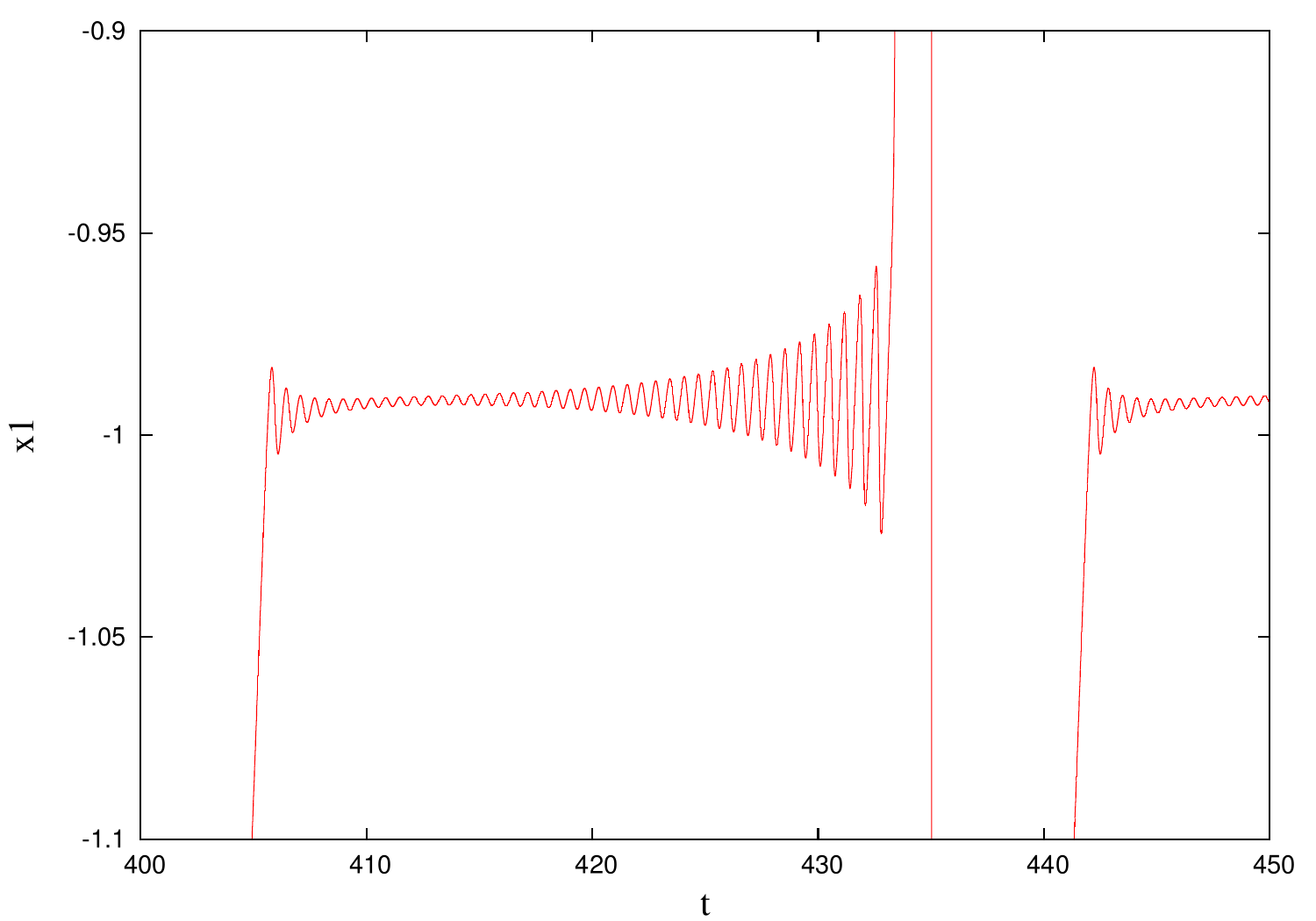}}
\subfigure[]{\includegraphics[scale=.3,angle=-90]{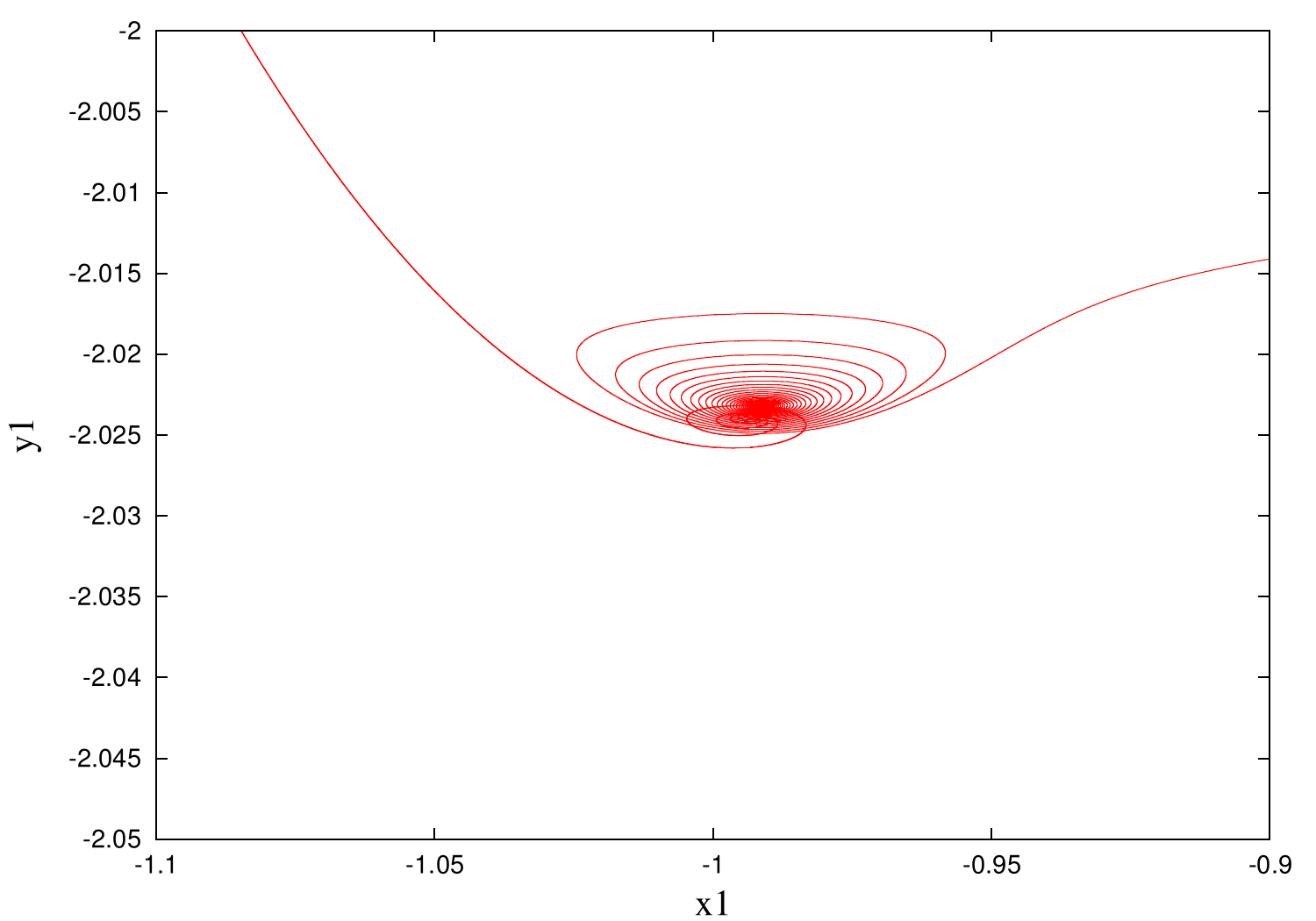}}\\
\caption{This figure corresponds to the simulation of system \eref{FHNdl} for a value of parameter $c_1=-0.952$. It is the same simulation that one represented in figure \ref{plc2Dfig}, but with a zoom  on the zone of small oscillations.  In panel (a), we show the time-evolution of $x_1$, for $t$ in the interval [400,450] and $x_1$ in the interval [-1.1,-0.9]  . In panel (b), we show the behavior of the system  in the $(x_1, y_1)$ phase-plane, for $x_1$ in the interval [-1.1,-0.9] and $y_1$ in the intervall [-2.05,-2]. We can see, the trajectory on the attractive manifold, then the small oscilations, resulting from the intersection of attractive and repulsive manifolds, and finally, the trajectory along the repulsive manifold before the evolution toward the fast direction.   }
\label{MMO2Dfig}
\end{center}
\end{figure}

\begin{figure}
\begin{center}
\includegraphics[scale=.5,angle=-90]{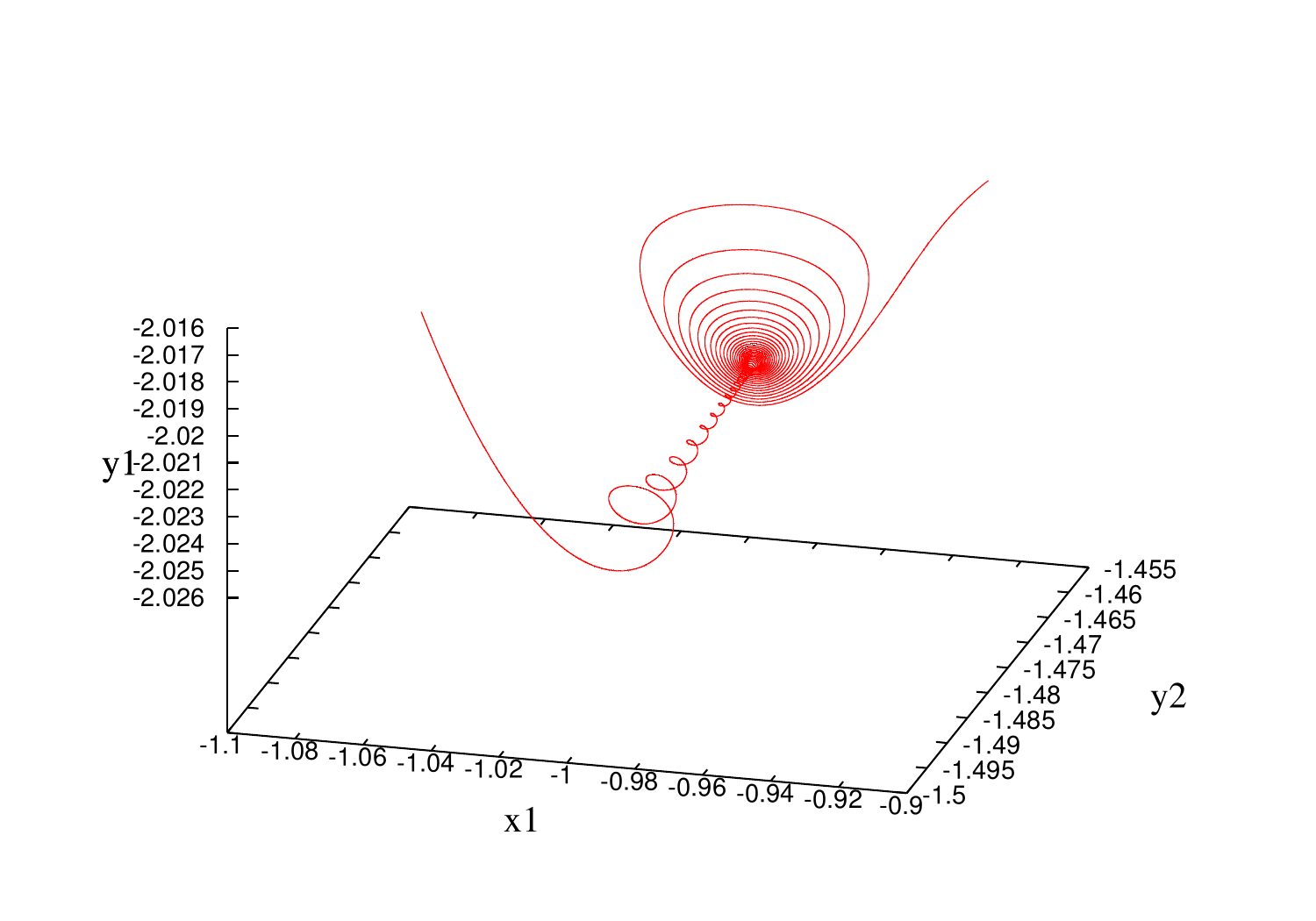}
\caption{This figure corresponds to the simulation of system \eref{FHNdl} for a value of parameter $c_1=-0.952$. It is the same simulation that the one represented in figure \ref{plc2Dfig}, but with a zoom  on the area of small oscillations.   We show the behavior of the system  in the $(x_1,y_2, y_1)$ phase-plane, for $x_1$ in the interval [-1.1,-0.9] and $y_1$ in the intervall [-2.026,-2.016]. Here we can easily distinguish the behavior described below. The small oscillations occur when the trajectories originating in the attracting manifold  in one side of the strong canard are trapped by the repelling one and have to return to the attracting manifold.  On the other side of the folded node, the trajectories follow the fast direction, see \cite{mb-mk-mw_06}. }
\label{MMO2Dbfig}
\end{center}
\end{figure}

\begin{figure}
\begin{center}
\includegraphics[scale=0.4, angle=-90]{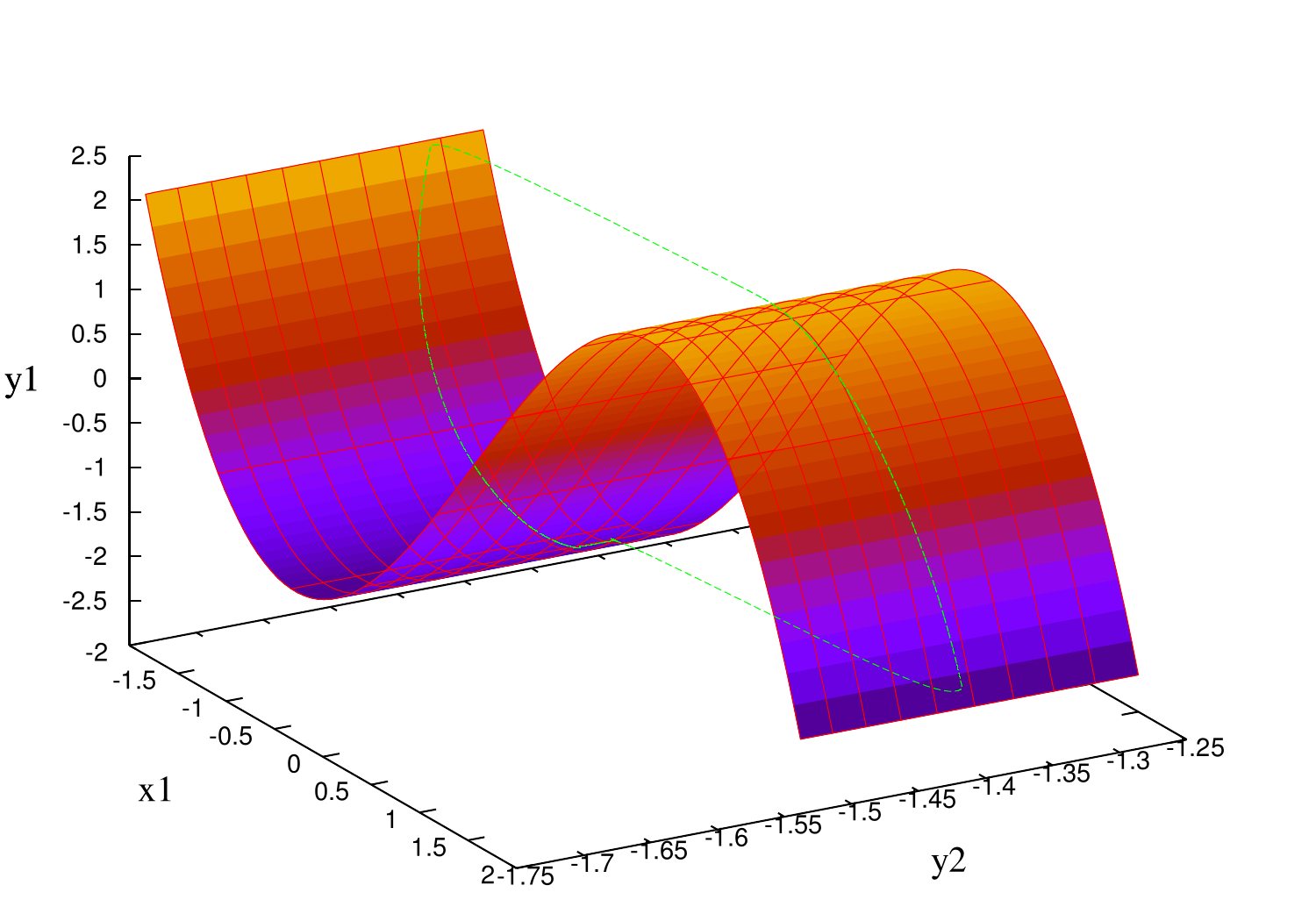}
\caption{This figure shows the simulation of system \eref{FHNdl} for parameter $c_1=-0.952$. We can see the critical manifold defined by equations: $(x_1,y_2=F(x_2)+\alpha_2(x_1-x_2),y_1=F(x_1)+\alpha_1(x_2-x_1)$ and  the behavior of the system in the $(x_1,y_2,y_1)$ phase-space (green curve). It shows the trajectory along the attractive critical manifold, then the evolution near the fold curve, the evolution along the fast direction and the global return macanism.}
\label{MMO3da}
\end{center}
\end{figure}

\begin{figure}
\begin{center}
\includegraphics[scale=0.3, angle=-90]{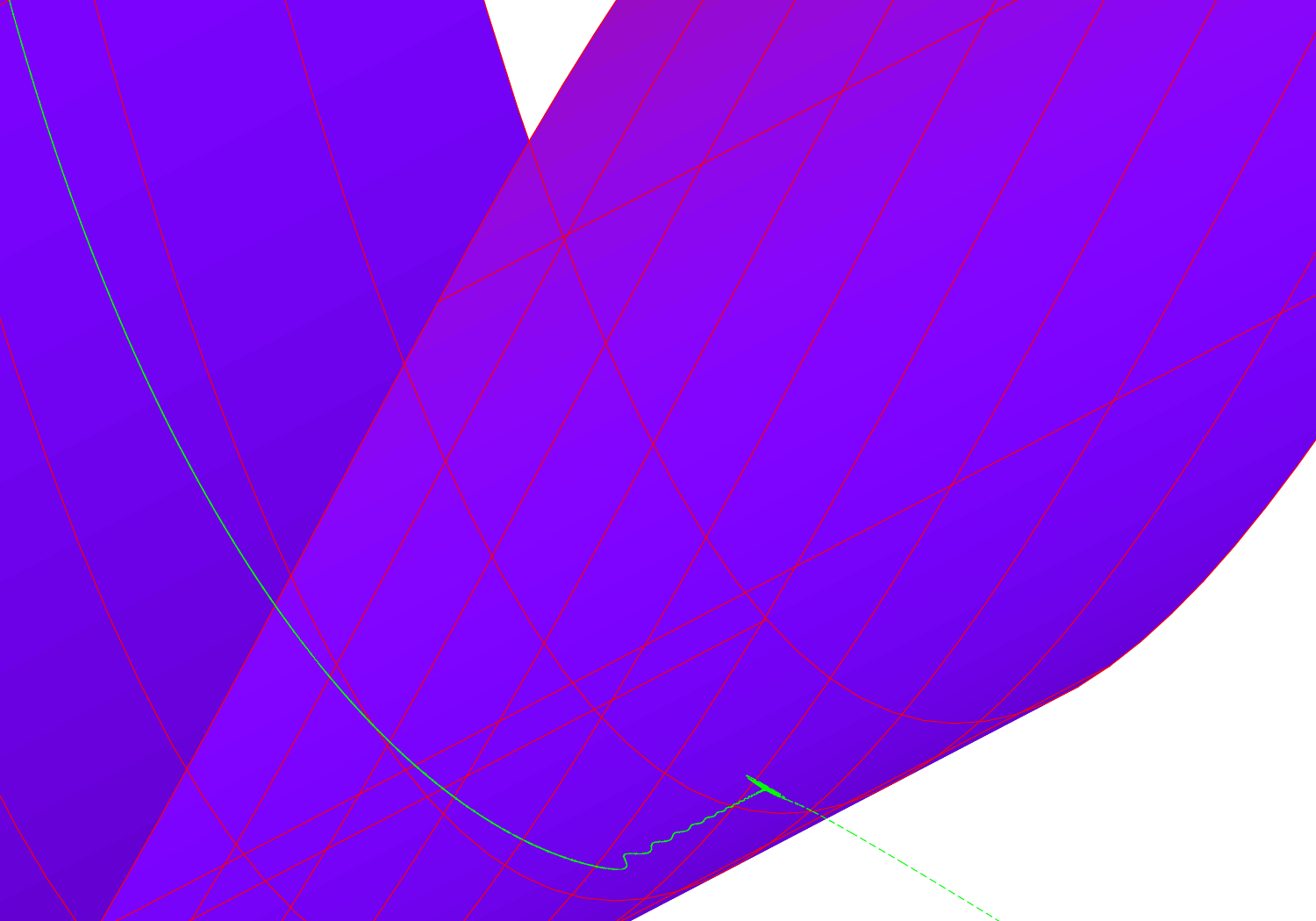}
\caption{This figure shows the simulation of system \eref{FHNdl} for parameter $c_1=-0.952$. We can see the critical manifold defined by equtions: $(x_1,y_2=F(x_2)+\alpha_2(x_1-x_2),y_1=F(x_1)+\alpha_1(x_2-x_1))$ and  the behavior of the system in the $(x_1,y_2,y_1)$ phase-space (green curve). It shows the trajectory along the attractive critical manifold, then the evolution near the fold curve and finally the evoluton along the fast direction.}
\label{MMO3db}
\end{center}
\end{figure}

\begin{figure}
\begin{center}
\includegraphics[scale=0.3, angle=-90]{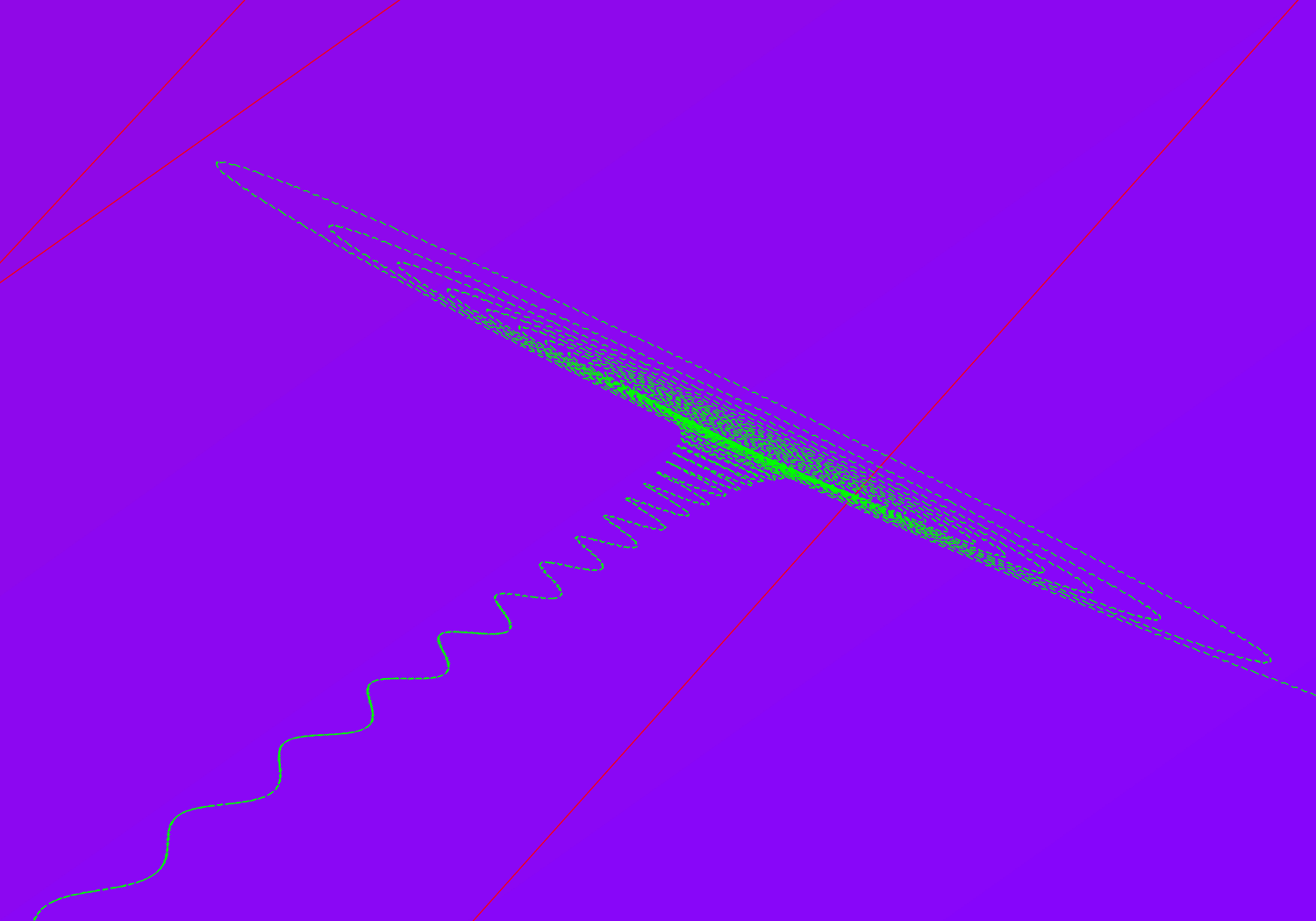}
\caption{This figure shows the simulation of system \eref{FHNdl} for parameter $c_1=-0.952$. We can see the critical manifold defined by equtions: $(x_1,y_2=F(x_2)+\alpha_2(x_1-x_2),y_1=F(x_1)+\alpha_1(x_2-x_1))$ and  the behavior of the system in the $(x_1,y_2,y_1)$ phase-space (green curve). It shows the small ocsilations occuring when the trajectory oscillates around the weak canard. The amplitute of oscillations when entering in the funnel region start to decrease before increasing until they cross the folded node and follow the fast direction.}
\label{MMO3dd}
\end{center}
\end{figure}

 Then, we vary the parameter $c_1$ in a small neighborhood around the value of $\bar{c}_1$. If $c_1<\bar{c}_1$, the system goes to a stationary state. As $c_1$ crosses the value $\bar{c}_1$ we can observe the appearance of MMOs. Below, we illustrate the simulation of system \eref{FHNdl} for two values of $c_1$. For $c_1=-0.953$, wich corresponds to a case where the system goes to the stationary state and for $c_1=-0.952$, wich corresponds to a case where we observe MMOs.  Furthermore, we can approximate, for the folded stationary point, the  eigenvalues of the Jacobian of \eref{Deseq1}-\eref{Deseq2}. For $c_1=-0.953$, we find:
 \[\lambda_1 \simeq 0.998, \ \lambda_2\simeq -0.02 .\]
 And ror $c_1=-0.952$, we find:
 \[\lambda_1 \simeq 1.004, \ \lambda_2\simeq 0.004.\]
 Following  \cite{mb-mk-mw_06}, \cite{ScholarWech:2007}, this gives an approximate theoretical number of small oscillations:
 \[s=[\frac{1}{2}+\frac{\lambda_1}{\lambda_2}]\simeq 135.\]
 In figure \ref{pointstable2Dfig}, we have illustrated the results of the simulation, for the parameter value $c_1=-0.953$. It shows, that in this case the system evolves asymptotically towards his stationnary state. The other figures illustrate the simulation results for $c_1=-0.952$. In Figure \ref{plc2Dfig} panel (a), we show the time evolution of the fast variable $x_1$, whereas in Figure \ref{plc2Dfig} panel (b) we show the behavior in the $(x_1,y_1)$ phase-plane. Here the MMOs appear but the small oscillations are hardly distinguishable because of their tiny amplitude in comparison with the relaxation oscillations. So, This, in figures \ref{MMO2Dfig} and \ref{MMO2Dbfig} we show a zoom of the previous illustration. In \ref{MMO2Dfig} panel (a), we show the time evolution of the fast variable $x_1$, whereas in Figure \ref{MMO2Dfig} panel (b), we show the behavior in the $(x_1,y_1)$ phase-plane. In panel (a), we can easily distinguish the small oscillations. In panel (b), we  can see the trajectory on the attractive manifold, then the small oscillations, resulting from the intersection of attractive and repulsive manifolds, and finally, the trajectory along the repulsive manifold before the evolution toward the fast direction. Figure \ref{MMO2Dbfig} is interesting as it clearly illustrates, in the $(x_1,y_2,y_1)$ phase plane, the MMOs in the case of  folded node singularity. The trajectories of system close to the singular funnel enter a region near the fold where they rotate around the weak canard: the trajectories originating in the attracting manifold  in one side of the strong canard are trapped by the repelling one and have to return to the attracting manifold, see \cite{mb-mk-mw_06}, for details. Finally, figures \ref{MMO3da},\ref{MMO3db} and \ref{MMO3dd}, also show the MMOs in the space $(x_1,y_2,y_1)$. In figure \ref{MMO3da}, we show the crititical manifold defined by $(x_1,y_2=F(x_2)+\alpha_2(x_1-x_2),y_1= F(x_1)+\alpha_1(x_2-x_1))$, and the trajectory of the system in the $(x_1,y_2,y_1)$ phase-space. We can see the small oscillations near the fold, the trajectory along the fast direction and the global return mechanism. The figure \ref{MMO3db} is the same as the figure \ref{MMO3da}, but with a zoom on area of the fold.  The same is done, in figure \ref{MMO3dd} but with a greater zoom in the area of small oscillations.

\section{An alternative approach to the analysis of coupled oscillator systems}
In this section we outline an alternative approach that was pointed out to us by the anonymous referee.
The modification consists of defining the reduced flow of \eref{FHNdl} by means of the fast variables
$x_1$ and $x_2$. Differentiating the equation $l(x_1,x_2,y_1,y_2)=0$ (see \eref{f} for the definition) we obtain 
\begin{center}
\begin{equation}
\label{doty}
\left ( \begin{array}{c} \dot y_1\\\dot y_2\end{array}\right )
=\left   (
       \begin{array}{cc}
      F'(x_1)-\alpha_1 & \alpha_1 \\
       \alpha_2& F'(x_2)-\alpha_2  \end{array}  \right  )  
      \left (
     \begin{array}{c} \dot{x}_1\\\dot{x}_2\end{array}
      \right )
\end{equation}
\end{center}
Hence, the reduced equation can be written in the form:
\begin{center}
\begin{eqnarray}
\label{eq-red2fast}
&\left (
    \begin{array}{cc}
      F'(x_1)-\alpha_1 & \alpha_1 \\
       \alpha_2& F'(x_2)-\alpha_2  \end{array}
      \right )  \left (
     \begin{array}{c} \dot{x}_1\\\dot{x}_2\end{array}\right )\\\nonumber&\qquad=\left (
     \begin{array}{c} x_1-c_1+\beta_1(F(x_2)-F(x_1)-(\alpha_1+\alpha_2)(x_2-x_1))\\
      x_2-c_2+\beta_2(F(x_2)-F(x_1)-(\alpha_1+\alpha_2)(x_2-x_1))
     \end{array}
      \right ).
\end{eqnarray}
\end{center}
We multiply both sides of \eref{eq-red2fast} by the cofactor matrix
\begin{equation}\label{eq-cofama}
\left (
     \begin{array}{cc}
      F'(x_2)-\alpha_2 & -\alpha_1 \\
       -\alpha_2& F'(x_1)-\alpha_1  \end{array}
      \right ),
\end{equation}
which yields
\begin{eqnarray}
\label{eq-red}
&((F'(x_1)-\alpha_1)(F'(x_2)-\alpha_2)-\alpha_1\alpha_2)  \left (
     \begin{array}{c} \dot{x}_1\\\dot{x}_2\end{array}\right )\\=&
      \left (
     \begin{array}{cc}
      F'(x_2)-\alpha_2 & -\alpha_1 \\
       -\alpha_2& F'(x_1)-\alpha_1  \end{array}
      \right )\left (
     \begin{array}{c} x_1-c_1+\beta_1(F(x_2)-F(x_1)-(\alpha_1+\alpha_2)(x_2-x_1))\\
      x_2-c_2+\beta_2(F(x_2)-F(x_1)-(\alpha_1+\alpha_2)(x_2-x_1))
     \end{array}
      \right ).\nonumber
\end{eqnarray}
Now we can desingularize by canceling the factor $(F'(x_1)-\alpha_1)(F'(x_2)-\alpha_2)-\alpha_1\alpha_2$,
which results in the following desingularized system
\begin{center}
\begin{eqnarray}
\label{eq-des2fast}
 & \left (
     \begin{array}{c} \dot{x}_1\\\dot{x}_2\end{array}\right )=\\&
      \left (
     \begin{array}{cc}
      F'(x_2)-\alpha_2 & -\alpha_1 \\
       -\alpha_2& F'(x_1)-\alpha_1  \end{array}
      \right )\left (
     \begin{array}{c} x_1-c_1+\beta_1(F(x_2)-F(x_1)-(\alpha_1+\alpha_2)(x_2-x_1))\\
      x_2-c_2+\beta_2(F(x_2)-F(x_1)-(\alpha_1+\alpha_2)(x_2-x_1))
     \end{array}
      \right ).\nonumber
\end{eqnarray}
\end{center}

\noindent This is equivalent to the desingularization that takes (7)-(8) to (9)-(10). 
In the context of \eref{eq-des2fast} folded singularities are equilibria on the fold line defined by $((F'(x_1)-\alpha_1)(F'(x_2)-\alpha_2)-\alpha_1\alpha_2)=0$.
Note that such equilibria exist robustly as  the cofactor matrix is singular along the fold line. 
To verify that we have a FSN II we have to find a transcritical bifurcation, with one equilibrium
crossing the fold line and the other one staying on it.  

This approach seems viable and may have computational advantages
over the one we have used. It is suited to the context of coupled oscillators, but relies
on the fact that there are at least as many fast variables as there are slow variables.
\section{Conclusion}
In this article, we have studied a system of two coupled slow-fast oscillators. We showed that coupling them leads to the occurrence of mixed mode oscillations. The main observation was that coupling the two systems have rise to an FSN II point, wich is known to imply the existence of small amplitude oscillations in the fold region. 
Our study has been the first attempt to understand canards and MMOs in systems of coupled oscillators. Our future goal is to extend our work to the context of large systems of oscillators. As coupled oscillators systems arise as discretizations of Reaction Diffusion equations, we hope to use the insights of this and future work to understand MMOs in the context of Reaction Diffusion equations. 
As numerous models in biology and neuroscience are constructed as networks of oscillators and Reaction Diffusion equations, this research has to be done  with deep interactions with applications.
\ack
We would like to thank: Région Haute Normandie  and
- FEDER (RISC project) for financial support.
\clearpage
\section*{References}
\bibliographystyle{plain}
\bibliography{Ambrosio}
\end{document}